\pdfoutput=1
\RequirePackage{ifpdf}
\ifpdf
\documentclass[pdftex]{sigma}
\else
\documentclass{sigma}
\fi

\usepackage{mathtools}

\usepackage[dvipsnames]{xcolor}

\let\ve\mathbf
\let\epsilon\varepsilon

\def\eqref#1{{\rm(\ref{#1})}}
\def\dif{\mathop{}\!\mathrm d}
\def\RR{{\mathbb R}}

\def\pd#1#2{\frac{\partial#1}{\partial#2}}

\def\P{{\mathcal P}}
\def\Q{{\mathcal Q}}

\def\U{{\mathcal U}}
\def\V{{\mathcal V}}

\def\RO{{\mathcal R}}
\def\Z{{\mathcal Z}}
\def\fpb#1#2{\fbox{\parbox{#2}{\small\centering #1}}}
\def\fpa#1#2{\fpb{$\begin{array}{r@{\,}l} #1 \end{array}$}{#2}}

\def\pmb#1{\hbox{\kern0.3pt#1\llap{#1\kern0.3pt}%
\kern0.3pt\llap{#1\kern0.3pt}\kern0.3pt\llap{#1\kern0.3pt}}}
\def\boldphi{\pmb{$\phi$}}

\numberwithin{equation}{section}

\newtheorem{Theorem}{Theorem}[section]
\newtheorem*{Theorem*}{Theorem}
\newtheorem{Corollary}[Theorem]{Corollary}
\newtheorem{Lemma}[Theorem]{Lemma}
\newtheorem{Proposition}[Theorem]{Proposition}

\theoremstyle{definition}
\newtheorem{Definition}[Theorem]{Definition}

\newtheorem{Example}[Theorem]{Example}
\newtheorem{Remark}[Theorem]{Remark}

\begin{document}

\allowdisplaybreaks

\newcommand{\arXivNumber}{2511.15558}

\renewcommand{\PaperNumber}{075}

\FirstPageHeading

\ArticleName{Voss Surfaces in Sine-Gordon Hierarchies}

\ShortArticleName{Voss Surfaces in Sine-Gordon Hierarchies}

\Author{Michal MARVAN}

\AuthorNameForHeading{M.~Marvan}

\Address{Mathematical Institute in Opava, Silesian University in Opava, \\ Na Rybn\'\i\v{c}ku 1, 746 01 Opava, Czech Republic}
\Email{\mail{michal.marvan@math.slu.cz}}

\ArticleDates{Received November 20, 2025, in final form July 23, 2026; Published online August 14, 2026}

\Abstract{We explore a~method, initiated by Guichard in 1890, which allows to generate sequences of Voss surfaces, starting from an arbitrarily chosen pseudospherical surface and a~seed solution of the Moutard equation, by means of two simple transformations. In this paper, we 1)~identify the Guichard transformations with the well-known recursion operator for symmetries of the sine-Gordon equation and its inverse; 2)~prove a~lemma which allows us to derive the length of Guichard's sequences from the invariance properties of the initial sine-Gordon solution; 3)~introduce an extended class of inverted operators, expanding the class of Voss surfaces obtainable by quadratures; 4)~clarify relevant aspects of Guthrie's formalism, paving the way for the future employment of the entire division algebra of recursion operators. A~number of Voss nets are presented explicitly.\looseness=-1}


\Classification{35A30; 37K10; 53A05; 58J70}

\section{Introduction}

In 1888, A. Voss studied conjugate nets formed by
geodesics~\cite{AV-1888}, uncovering their basic properties and
their relation to pseudospherical surfaces, thereby also
to the sine-Gordon equation
\begin{gather}
\label{eq:sG}
\phi_{xy} = \sin\phi.
\end{gather}
Independently, C.~Guichard~\cite{CG-1890b} investigated the
same class in the context of congruences, being the first to notice
the relation to the Moutard equation
\begin{equation}
\label{eq:VsG}
\Phi_{xy} = \Phi \cos \phi
\end{equation}
and its transformations.
Accounts of the first period of investigation were given by
Gambier~\cite{BG-1931}
and to some extent also by Salkowski~\cite{ES-1924}.
The problem of finding explicit solutions was addressed by
many authors, notably
Eisenhart~\cite{LPE-1906,LPE-1909,LPE-1914,LPE-1915,LPE-1917}.
Razzaboni~\cite{AR-1889} sought minimal Voss surfaces.
The most familiar explicit example is the right helicoid~\cite{AV-1888},
carrying infinitely many Voss nets~\cite{SF-1939,BG-1926}.
For the Voss surfaces of revolution see Tachauer~\cite{AT-1903}
and Gambier~\cite{BG-1927,BG-1931},
for explicit Voss surfaces associated with the Enneper surfaces see
Diller~\cite{JBD-1913}.

Voss surfaces are also amenable to modern methods of integrable
systems.
Babich~\cite{Bab-1996} obtained Voss surfaces from the
$\psi$-function of the Lax pair for the sine-Gordon equation.
In a~different way, Lax integrability was established by Lin and
Conte~\cite{RL-RC-2018}, who also revisited particular
surfaces considered by Gambier.

Voss nets have seen a~recent resurrection
in connection with geometric modelling.
In large-span free-form architecture, discrete Voss nets
(Sauer and Graf~\cite{RS-HG-1931}, Sauer~\cite{RS-1970},
Doliwa~\cite{AD-1999}, Izmestiev et al.\ \cite{Izm})
model elastically bent lath constructions with planar cladding
(Montagne et al.\ \cite{MDTFB-2020,MDTFB-2022},
Kilian et al.\ \cite{KNRRS-2024} and references therein).
Smooth Voss nets are smooth limits of
discrete Voss nets, revealing shapes achievable with this
construction technology.
In this vein, Rasoulzadeh~\cite{AR-2025} investigated
deployable Voss nets.

In this paper, we draw on linking the classical results
with the modern theory of integrable systems.
Our starting point is the following classical construction,
summarised in Section~\ref{sect:Voss}.
In short,
any solution $\phi$ to the sine-Gordon equation~\eqref{eq:sG} determines a~pseudospherical surface via the Gauss--Weingarten equations;
any solution $\Phi$ to the Moutard equation~\eqref{eq:VsG} then
determines a~Voss surface algebraically,
although possibly a~degenerate one.
This extends to linear combinations of solutions $\Phi$ by linearity.

In Section~\ref{sect:sG2Voss}, we observe that Moutard's equation
essentially matches the sine-Gordon symmetry condition.
This is the way symmetries of the sine-Gordon equation give rise to Voss
surfaces.
We work through three simple explicit examples, summarised in
Table~\ref{tab:ex}.

In Section~\ref{sect:sG:RO}, we recall Olver's~\cite{PJO-1977}
recursion operator $\RO$ for symmetries of the sine-Gordon equation.
By definition, this is an integro-differential operator that sends
symmetries to symmetries.
To resolve multiple issues inherent in the integro-differential framework,
Guthrie~\cite{GAG-1993,GAG-1994} reinterpreted recursion operators
essentially as B\"acklund auto-transformations for the linearised
equation~\cite{MM-alro}.\footnote{This point of view was actually pioneered in a~series of papers by
Papachristou, e.g.,~\cite{CJP-1991}.
It is quite explicit in recent works by Habibullin et al., e.g.,
\cite{H-Kh-P-2016}.}

Incidentally, the Guthrie form of $\RO$ coincides with a~transformation discovered by Guichard in 1890,
the same being true for the inverse $\RO^{-1}$.
Since the treatment of $\RO^{-1}$ appears to be somewhat
problematic in recent literature
(see footnote~\ref{foo:Dinv} below),
we pay special attention to understanding of $\RO^{-1}$.

When iterated, both $\RO$ an $\RO^{-1}$ generate two infinite families of
symmetries, namely, the pmKdV and the Khor'kova hierarchies
in two versions, differing by interchanging $x$ and $y$.
These are discussed in Section~\ref{sect:pos}, along with some examples
of the Guichard sequences of Voss nets, obtained from one and the same
pseudospherical surface with the help of the recursion operator.
Examples show that the sequences stop after a~few steps,
as observed already by Guichard.\looseness=-1

The stopping problem (we call it the Guichard problem) is addressed in
Section~\ref{sect:GP}.
With our Lemma~\ref{lemma1}, we link the number of independent
nontrivial Voss surfaces in the Guichard sequence to the symmetry
invariance properties of the initial sine-Gordon solution.

While Lemma~\ref{lemma1} limits the number of Voss
surfaces obtainable with the help of $\RO$ and $\RO^{-1}$,
it imposes no restrictions on results obtained with the help of the
general inversion $(\RO + \lambda\,\mathrm{Id})^{-1}$.
These are examined in Section~\ref{sect:neg} after reduction to
$(\RO + \mathrm{Id})^{-1}$.
The main result is a~formula for
$(\RO + \mathrm{Id})^{-1} \Phi$
 given in
Proposition~\ref{prop:sG:iro}.
We explicitly obtain the first member of what we call the negative
$\RO + \mathrm{Id}$ hierarchy.

In Section~\ref{sect:iRO:Voss}, we demonstrate that the new hierarchy
does indeed take us outside of the Guichard sequence.
What is more, we show that Proposition~\ref{prop:sG:iro} opens another
path from a~pseudospherical surface to a~Voss net. A~chance observation of this new route was the initial impetus for writing
the present paper.

The paper is closed with a~digression on the algebra of recursion operators.
We prove that every element of the division algebra generated by $\RO$
yields a~true recursion operator for the sine-Gordon equation.
This prepares the ground for the future use of the entire division algebra.

Smoothness is assumed throughout the paper.

\section{Essentials about Voss nets and Voss surfaces}
\label{sect:Voss}

In this section, we briefly review known geometric facts that link
Voss surfaces to pseudospherical surfaces and solutions of the
Moutard equation.

According to Dini~\cite{UD-1870}, a~surface is pseudospherical (i.e., has a~constant negative Gaussian curvature)
if and only if it admits a~parameterisation such that the fundamental forms are
\begin{gather}
\mathrm{I} = \dif x^2 + 2 \cos \phi \dif x \dif y + \dif y^2,
\qquad
\mathrm{II} = 2 \sin \phi \dif x \dif y,
\nonumber
\\
\mathrm{III} = \dif x^2 - 2 \cos \phi \dif x \dif y + \dif y^2,
\label{eq:sG:ff}
\end{gather}
where $\phi$ is the {\it intersection angle}, also known as the
{\it Chebyshev angle}.
The parameterisation~\eqref{eq:sG:ff} will be called
the {\it Dini parameterisation}.
The form of $\mathrm{I}$ (resp.\ $\mathrm{II}$) implies that the Dini
parameterisation is {\it Chebyshev} (resp.\ {\it asymptotic}) on the
surface.
The form of $\mathrm{III}$ conveys that Dini's parameterisation is
Chebyshev on the Gaussian image.

The Gauss--Weingarten system corresponding to the
fundamental forms~\eqref{eq:sG:ff} is
\begin{gather}
\ve r_{xx} = \frac{\phi_x}{\tan \phi} \ve r_x - \frac{\phi_x}{\sin \phi} \ve r_y,
\qquad
\ve r_{xy} = \sin \phi \ve n,
\qquad
\ve r_{yy} = -\frac{\phi_y}{\sin \phi} \ve r_x + \frac{\phi_y}{\tan \phi} \ve r_y,
\nonumber
\\
\ve n_x = \frac{1}{\tan \phi} \ve r_x - \frac{1}{\sin \phi} \ve r_y,
\qquad
\ve n_y = \frac{1}{\tan \phi} \ve r_y - \frac{1}{\sin \phi} \ve r_x,
\label{eq:sG:GW}
\end{gather}
while the Gauss--Mainardi--Codazzi system reduces to the
sine-Gordon equation~\eqref{eq:sG}.
Thus, the sine-Gordon equation determines both pseudospherical
surfaces and Chebyshev nets on the sphere.

\begin{Definition}[\cite{AV-1888}]
{\it Voss nets} are conjugate nets formed by geodesics.
{\it Voss surfaces} are surfaces capable of carrying a~Voss net.
\end{Definition}

Recall that conjugate nets are defined by the diagonality of
the second fundamental form, that is, $\mathrm{II}_{12} = 0$.
As already observed by Voss~\cite{AV-1888},
the Gaussian image of a~Voss net is a~Chebyshev net on the
Gauss sphere, which corresponds to a~solution of the sine-Gordon equation.
The problem of reconstruction of the Voss net
from its Gaussian image can be reduced to a~second-order linear PDE,
available in several distinct forms.\footnote{For instance,
Diller uses the linear equation
$M_{xy} \sin\phi + \phi_x M_y + \phi_y M_x = 0$
and finds solutions $M$ such that~$M_x$ and $M_y$ are functions of
$\phi$.}
One of them is the Moutard equation~\eqref{eq:VsG}, which lies
at the roots of the classical relationship between infinitesimal
deformations of pseudospherical surfaces and Voss nets in tangential
coordinates~\cite{Bia-II,LPE-1909,Izm}.
This piece of classical theory is summarised in
Proposition~\ref{VsG->Voss}, which is a~variation of~\cite[Theorem~2.25]{Izm}.
Explanatory references are given in Remark~\ref{rem:isodef}.

\begin{Proposition}
\label{VsG->Voss}
Let $\ve r(x,y)$ be a~pseudospherical surface in Dini's parameterisation,
let $\phi$ denote the Chebyshev angle.
Let $\Phi$ be a~solution of the Moutard equation
$\Phi_{xy} = \Phi \cos\phi$.
\begin{itemize}
\itemsep=0pt
\item[$(i)$]
The net
\begin{equation}
\label{eq:Voss:Phi}
\ve q = \Phi \ve n
 - \frac{\Phi_{y} \ve r_{x} + \Phi_{x} \ve r_{y}}{\sin\phi},
\end{equation}
if non-singular, is a~Voss net, possessing the same normal vector $\ve n$
and satisfying $\ve q \cdot \ve n = \Phi$.
Its fundamental forms are
\begin{equation}
\label{eq:Voss:FF}
\mathrm{I} = \frac{X^2 \dif x^2 + 2 X Y \cos\phi \dif x \dif y +Y^2 \dif y^2}
 {\sin\phi^2},
\qquad
\mathrm{II} = -X \dif x^2 - Y \dif y^2,
\end{equation}
where
\begin{equation}
\label{eq:Phi:XY}
X = \Phi_{xx} + \Phi
 - \left(\frac{\Phi_x}{\tan\phi} + \frac{\Phi_y}{\sin\phi}\right) \phi_x,
\qquad
Y = \Phi_{yy} + \Phi
 - \left(\frac{\Phi_y}{\tan\phi} + \frac{\Phi_x}{\sin \phi}\right) \phi_y.
\end{equation}
The Gauss and mean curvatures are
\[
K = \frac{\sin^2\phi}{X Y}, \qquad
H = -\frac 1 2\left(\frac 1X + \frac 1Y\right),
\]
respectively.

\item[$(ii)$]
The singularity condition $\det\mathrm{I} = 0$ of the Voss
net~\eqref{eq:Voss:Phi}
is
\begin{equation}
\label{eq:Phi:sing}
X Y = 0.
\end{equation}

\item[$(iii)$]
If one of functions $X$, $Y$ vanishes everywhere,
then so does the other.
This happens if and only if the associated net degenerates to a~point.

\item[$(iv)$]
Every non-singular Voss net is obtainable in the way described in
statement {\rm(i)}
from a~solution $\phi$ of the sine-Gordon equation
and a~solution $\Phi$ of the Moutard equation.
\end{itemize}
\end{Proposition}

\begin{proof}
(i) The fundamental forms I and II of the net $\ve q$ given
by equation~\eqref{eq:Voss:Phi} are obtained by straightforward computation.
Since $\mathrm{II}_{12} = 0$, the parametric net is conjugate.
Computation of geodesic curvatures shows that the net is also
geodesic. The curvatures $K$, $H$ are also straightforward to compute.

(ii) Clearly,
$\det\mathrm{I} = X Y/{\sin^2\phi}$,
which proves the second statement.

(iii) It is easy to check that
\begin{equation}
\label{eq:Voss:XYsys}
X_y = -\frac{\phi_x}{\sin\phi} Y,
\qquad
Y_x = -\frac{\phi_y}{\sin\phi} X
\end{equation}
in consequence of the sine-Gordon and Moutard equations.
Hence, if $X = 0$, then $Y = 0$, and vice versa.
Moreover,
\begin{equation}
\label{eq:Voss:Phi:TD}
\ve q_x = -\frac{X}{\sin\phi} \ve r_y,
\qquad
\ve q_y = -\frac{Y}{\sin\phi} \ve r_x,
\end{equation}
which finishes the proof of the third statement.

(iv)
Consider a~conjugate net with the first fundamental form I given by
formula~\eqref{eq:Voss:FF} and the second fundamental form
$\mathrm{II} = \xi X \dif x^2 + \eta Y \dif y^2$, where $X$, $Y$, $\xi$, $\eta$ are
arbitrary functions of~$x$,~$y$;
this is generic for conjugate nets.
The condition that the isoparametric curves are geodesics is
equivalent to system~\eqref{eq:Voss:XYsys}.
Assuming~\eqref{eq:Voss:XYsys}, the Gauss--Mainardi--Codazzi system
reduces to~${
\xi_{y} = 0}$,
$
\eta_{x} = 0$,
$
\phi_{xy} = \xi \eta \sin \phi$.
Hence, $\xi$ and $\eta$ are functions of $x$ and $y$, respectively.
As such, they can be removed by separate transformations of
$x$ and $y$, preserving the net.
We are left with the sine-Gordon equation for $\phi$ and
system~\eqref{eq:Voss:XYsys} for $X$, $Y$.
These imply compatibility of the system
\begin{gather*}
\Phi_{xx} + \Phi
 - \left(\frac{\Phi_x}{\tan\phi} + \frac{\Phi_y}{\sin\phi}\right) \phi_x = X,
\qquad
\Phi_{xy} - \Phi \cos\phi = 0,
\\
\Phi_{yy} + \Phi
 - \left(\frac{\Phi_y}{\tan\phi} + \frac{\Phi_x}{\sin \phi}\right) \phi_y = Y
\end{gather*}
composed of the Moutard equation and equations~\eqref{eq:Phi:XY}
solved for $\Phi_{xx}$ and $\Phi_{yy}$.
Then the net $\ve q$ defined by formula~\eqref{eq:Voss:Phi} has the
same normal vector $\ve n$, the same support function $\Phi$, and
the same fundamental forms~\eqref{eq:Voss:FF}.
Therefore, $\phi$ and $\Phi$ are the solutions sought.
\end{proof}

The Voss net~\eqref{eq:Voss:Phi} is said to be {\it associated} with the
solution~$\Phi$ and the initial pseudospherical surface $\ve r$.\footnote{This is consistent with the well-known general definition of
associate nets \cite[Section~225]{Bia-II}, \cite[Section~155]{LPE-1909}
via the relation
$\mathrm{II}_{11} \mathrm{II}_{22}' - 2 \mathrm{II}_{12} \mathrm{II}_{12}' + \mathrm{II}_{22} \mathrm{II}_{11}' = 0$.}

\begin{Remark} \label{rem:isodef}
We refer the reader to Bianchi~\cite{Bia-I,Bia-II} or
Eisenhart~\cite{LPE-1909} or Izme\-stiev~\cite{Izm}
for a~broader geometric perspective.
In particular, the geometric meaning contained in
Proposition~\ref{VsG->Voss} can be explained as follows.
\begin{itemize}\itemsep=0pt
\item[(1)] According to \cite[Section~81]{Bia-I} or \cite[Section~67]{LPE-1909},
the {\it support function} of a~surface $\ve q$ is the
distance $\ve q \cdot \ve n$ of the tangent planes from the origin;
the normal vector $\ve n$ and the support function
$w = \ve q \cdot \ve n$ are referred to as the {\it tangential coordinates}
of a~surface.

\item[(2)] For a~pseudospherical surface $\ve r$ in Dini's parameterisation,
solutions $\Phi$ of the Moutard equation can be identified with
the Weingarten characteristic functions of infinitesimal isometric
deformations, see \cite[Chapter~XV, Section~226, case~$1^\circ$]{Bia-II}
or~\cite[Chapter~XI]{LPE-1909}
or~\hbox{\cite[Section~2.5]{Izm}}.
According to Bianchi \hbox{\cite[Section~225]{Bia-II}}, the associate net
$\ve q$ can be obtained by taking $\Phi$ for the support function $w$;
it is then known to possess the Voss property.

\item[(3)] Once $\Phi$ and $\ve n$ are known,
the position vector $\ve q$ can be obtained by solving the linear system
\begin{equation}
\label{eq:q:sys}
\ve q \cdot \ve n = \Phi, \qquad
\ve q \cdot \ve n_x = \Phi_x, \qquad
\ve q \cdot \ve n_y = \Phi_y
\end{equation}
(see \cite[Section~81, equations~(a) and~(b)]{Bia-I},
\cite[Section~67, equations~(29) and~(30)]{LPE-1909}).
The explicit expression for the solution of system~\eqref{eq:q:sys} is
\begin{equation}
\label{eq:q:nw}
\ve q = \mathrm{III}^{ij} \Phi_{,i} \ve n_{,j} + \Phi \ve n,
\end{equation}
where $\mathrm{III}^{ij}$ denotes the inverse of the Gauss sphere's metric
$\mathrm{III}_{ij}$,
cf.\ \cite[Chapter~V, Section~81, equation~(34)]{Bia-I} or
\cite[Section~67, equation~(32)]{LPE-1909}.

\item[(4)] Formula~\eqref{eq:q:nw},
which yields $\ve q$ in terms of $\Phi$ and $\ve n$,
becomes \hbox{\cite[equation~(24)]{Izm}},
which can be reduced to our formula~\eqref{eq:Voss:Phi} by expressing
$\ve n_x$ and $\ve n_y$ in terms of $\ve r_x$ and $\ve r_y$ via the
Mainardi--Codazzi equations.
This slight simplification makes the only difference between
Proposition~\ref{VsG->Voss} and \hbox{\cite[Theorem~2.25]{Izm}}.
\end{itemize}

The overall picture can be conveyed by the pullback diagram
$$
\def\fpb#1#2{\fbox{\parbox{#2}{\small\centering #1}}}
\begin{array}{ccccc}
\fpb{pseudospherical surface and its deformation}{4.2cm}
 & \xleftrightarrow{\text{characteristic $\leftrightarrow$ support}} & \fpb{Voss net}{1.2cm}
\\
\llap{\scriptsize }\downarrow
 & \hbox{\scriptsize pullback}
 & \downarrow \rlap{\scriptsize Gauss map}
\\
\fpb{pseudospherical surface}{3cm}
 & \xleftrightarrow[\text{Gauss--Lelieuvre}]{}
 & \fpb{Chebyshev net on the sphere}{2.7cm}
\end{array}
$$
By the deformation we mean the infinitesimal isometric deformation
\cite[Chapter~XV]{Bia-II}
or~\cite[Chapter~XI]{LPE-1909}.
\end{Remark}

\begin{Corollary}
\label{cor:Voss:PhiPhi'}
Given a~pseudospherical surface of Chebyshev angle~$\phi$,
let $\Phi$, $\Phi'$ be two solutions of the Moutard
equation~\eqref{eq:VsG} and
$X$, $X'$ $($resp.\ $Y$, $Y')$ be the corresponding functions~\eqref{eq:Phi:XY}.
Then the associated Voss nets
\[
\ve q = \Phi \ve n
 - \frac{\Phi_{y} \ve r_{x} + \Phi_{x} \ve r_{y}}{\sin\phi},
\qquad
\ve q' = \Phi' \ve n
 - \frac{\Phi'_{y} \ve r_{x} + \Phi'_{x} \ve r_{y}}{\sin\phi},
\]
differ by a~translation if and only if any of the equivalent conditions
$X = X'$ or $Y = Y'$ holds.
\end{Corollary}

\begin{proof}
By Proposition~\ref{VsG->Voss}\,(iii) and linearity.
\end{proof}

Given a~solution $\phi$ of the sine-Gordon equation, a~number of exact
solutions $\Phi$ of the Moutard equation are known.
Many have been known since the nineteenth century,
see Guichard~\cite{CG-1890b,CG-1890a}.
They will be reviewed in the next section.
However, not every pair $\phi$, $\Phi$ produces a~nontrivial Voss net
because of Proposition~\ref{VsG->Voss}\,(ii)
and Corollary~\ref{cor:Voss:PhiPhi'}.

\begin{Example}
\label{ex:nc}
For all pseudospherical surfaces $\ve r$,
the tangential coordinates $\Phi = \ve c \cdot \ve n$,
where $\ve c =$ const, and
$\Phi = w = \ve r \cdot \ve n$ are solutions of the Moutard
equation~\eqref{eq:VsG},
see~\cite[p.\ 234]{CG-1890b}.
Unfortunately, both satisfy the condition~\eqref{eq:Phi:sing}
and thus yield degenerate Voss nets.
According to Corollary~\ref{cor:Voss:PhiPhi'}, these are not
useful even in linear combinations with other solutions.
For another derivation see Proposition~\ref{prop:R+I:ideal} below.
\end{Example}

\section{Voss nets from symmetries}

\label{sect:sG2Voss}

In this section, we exploit the link between the Moutard equation and
symmetries of the sine-Gordon equation.
We do so because the sine-Gordon theory is
widely known, sufficiently developed, and well suited to our needs.
For a~background in the geometric theory of
symmetries, recursion operators and coverings one can consult
the paper~\cite{ISK-AMV-1989} or books~\cite{B-V-V,JK-AV-RV-2017}.

All the mixed derivatives
$\phi_{xy},\phi_{xxy},\phi_{xyy},\dots$
can be expressed in terms of the pure derivatives
$\phi,\phi_x,\phi_{xx},\phi_{xxx},\dots,\phi_y,\phi_{yy},\phi_{yyy},\dots$
via the sine-Gordon equation~\eqref{eq:sG};
the equation then determines a~manifold,
naturally coordinatised by
$x,y,\phi,\phi_x,\phi_{xx},\phi_{xxx},\dots,
 \phi_y,\phi_{yy},\phi_{yyy},\dots$ and
equipped with the {\it total derivatives,
that is, the vector fields
\[
D_x = \pd{}{x} + \sum_\mu \phi_{x\mu} \pd{}{\phi_\mu},
\qquad
D_y = \pd{}{y} + \sum_\mu \phi_{y\mu} \pd{}{\phi_\mu},
\]
where $\mu = x^m y^n$ denotes an arbitrary monomial,%
\footnote{The so-called {\it Janet monomials} are an obvious alternative
to the symmetric multiindices used in~\cite{B-V-V,JK-AV-RV-2017}.}
$\phi_{xy\mu} = D_\mu \sin\phi$, and
\[
D_{x^m y^n} = D_x^m \circ D_y^n = D_y^n \circ D_x^m.
\]
Assuming, as we do, that functions
$\boldphi(x,y,\phi,\phi_x,\phi_{xx},\dots,\phi_y,\phi_{yy},\dots)$
depend only on a~finite number of arguments,
all sums over $\mu$ in $D_x \boldphi$, $D_y \boldphi$ are finite.}

 A~function
$\boldphi(x,y,\phi,\phi_x,\phi_{xx},\dots,\phi_y,\phi_{yy},\dots)$
is called a~{\it symmetry}
(more precisely, a~{\it local symmetry}) of the sine-Gordon equation, if
\begin{equation}
\label{eq:VsG:D}
D_{xy} \boldphi = \boldphi \cos \phi.
\end{equation}
Equation~\eqref{eq:VsG:D} is called the
{\it linearised sine-Gordon equation in total derivatives}.
Examples of symmetries are provided by the members of the pmKdV
hierarchy~\eqref{eq:hierarchy:pmKdV}.

Denoting by $t$ an auxiliary variable different from $x$, $y$ (a ``time''),
the {\it flow} of a~symmetry $\boldphi$ is the evolution equation
\begin{equation}
\label{eq:flow}
\phi_t = \boldphi(x,y,\phi,\phi_x,\phi_{xx},\dots,\phi_y,\phi_{yy},\dots).
\end{equation}
Equation~\eqref{eq:VsG:D} is the compatibility condition between
the flow~\eqref{eq:flow} and the sine-Gordon equation~\eqref{eq:sG}.
Indeed, the standard cross-differentiation test gives
$D_{xy}\boldphi = \phi_{txy} = \phi_{xyt} = D_t(\sin\phi)\allowbreak = \phi_t \cos\phi
 = \boldphi \cos\phi$.

 A~solution $\phi$ of the sine-Gordon equation is said to be
{\it invariant} with respect to a~symmetry~$\boldphi$~if
$
\boldphi(x,y,\phi,\phi_x,\phi_{xx},\dots,\phi_y,\phi_{yy},\dots) = 0$.
Thus, the invariant solutions are just the stationary solutions of the
associated flows.

By {\it classical Lie symmetries} we mean symmetries
that depend on $x$, $y$, $\phi$ and linearly on~$\phi_x$,~$\phi_y$ at most.
For the sine-Gordon equation, they are listed in Table~\ref{tab:sG:sym}
along with their basic properties.
The scaling parameter $\gamma$ is also known as the
B\"acklund parameter.

\begin{table}
$$
\def\arraystretch{1}
\begin{tabular}{l|c|c|c}
symmetry $\boldphi$ & $\phi_x$ & $\phi_y$ & $x \phi_x - y \phi_y$
\\
name & \multicolumn2{c|}{translation} & scaling
\\
$1$-parametric group & $x \mapsto x + \alpha$ & $y \to y + \beta$ &
$x \mapsto \mathrm{e}^\gamma x$, $y \mapsto \mathrm{e}^{-\gamma} y$
\\
asymptotic parameterisation & \multicolumn2{c|}{preserved} & preserved
\\
Chebyshev parameterisation & \multicolumn2{c|}{preserved} & not preserved
\\
pseudospherical surface & \multicolumn2{c|}{preserved} & not preserved
\end{tabular}
$$
\caption{Classical Lie symmetries of the sine-Gordon equation.}
\label{tab:sG:sym}
\end{table}

\begin{Remark}\label{rem:sG:inv}
As is well known~\cite{P-M-2007},
there are two types of sine-Gordon solutions invariant
with respect to the classical symmetries.
They are
\begin{itemize}\itemsep=0pt
\item[(1)] The {\it travelling-wave solutions}, which are
solutions satisfying the invariant con\-straint $\phi_y = m \phi_x$,
where $m = {\rm const} \ne 0$.
They are of the form $\phi = \phi(x + m y)$.
Generic solutions of this kind are periodic and expressible in terms
of the Jacobi elliptic functions,
see \cite[Section~2.5.1]{P-M-2007} for details.
For the corresponding pseudospherical surfaces see~\cite{Z-1994}
and~\cite[Section~3.7.1]{P-M-2007}.
For gypsum models visit the corresponding section of the
G\"ottingen Collection of Mathematical Models and
Instruments~\cite[Section~E.V]{sammlung}.

The only travelling-wave solutions expressible in terms of elementary
functions are the kinks
\begin{equation}
\label{eq:Dini:phi}
4 \arctan\bigl(\mathrm{e}^{\mathrm{e}^\gamma x + \mathrm{e}^{-\gamma} y}\bigr)
\end{equation}
of velocity $-\mathrm{e}^{2\gamma}$, for arbitrary $\gamma \in \RR$.
They correspond to the {\it Dini helicoids}
\begin{gather}
\ve r = \frac1{\cosh\gamma}
\left[\frac{\cos(x - y)}{\cosh(\mathrm{e}^\gamma x + \mathrm{e}^{-\gamma} y)},
 \frac{\sin(x - y)}{\cosh(\mathrm{e}^\gamma x + \mathrm{e}^{-\gamma} y)},\nonumber
 -{\tanh}(\mathrm{e}^\gamma x + \mathrm{e}^{-\gamma} y)\right]
\\
\phantom{\ve r =}{} +
[0, 0, x + y].\label{eq:Dini:r}
\end{gather}
The corresponding spherical Chebyshev net is
\begin{gather*}
\ve n = -\frac{\tanh(\mathrm{e}^\gamma x + \mathrm{e}^{-\gamma} y)}{\cosh \gamma}
 \left[\cos(x - y), \sin(x - y),
 \frac{1}{\sinh(\mathrm{e}^\gamma x + \mathrm{e}^{-\gamma} y)}\right]
\\ \phantom{\ve n =}{} +
\tanh \gamma \left[-\sin(x - y), \cos(x - y), 0\right].
\end{gather*}

For $\gamma = 0$, the kink is
\begin{equation}
\label{eq:PS:phi}
\phi = 4\arctan \mathrm{e}^{x + y},
\end{equation}
while the Dini helicoid collapses to the Beltrami pseudosphere
\begin{equation}
\label{eq:PS:r}
\ve r =
\left[\frac{\cos(x - y)}{\cosh(x + y)}, \frac{\sin(x - y)}{\cosh(x + y)},
 x + y - \tanh(x + y)\right],
\end{equation}
and the spherical Chebyshev net to
\begin{equation}
\label{eq:PS:n}
\ve n =
-{\tanh}(x + y)
 \left[\cos(x - y), \sin(x - y), \frac{1}{\sinh(x + y)}\right].
\end{equation}

\item[(2)] The {\it scaling-invariant solutions}, that is,
the solutions satisfying the invariant constraint $x \phi_x - y \phi_y = 0$.
They are of the form $\phi = \phi(x y)$ and
expressible in terms of the third Painlev\'e transcendent.
For details see \cite[Section~2.5.2]{P-M-2007}.
The corresponding surfaces are known as the Amsler surfaces, see
\cite[Section~3.7.2]{P-M-2007} and references therein.
\end{itemize}
\end{Remark}

Obviously, solutions $\boldphi$ of
equation~\eqref{eq:VsG:D} can be viewed as differential operators
 (generally nonlinear)
that map solutions $\phi$ of the sine-Gordon equation~\eqref{eq:sG}
to solutions
\[
\Phi = \boldphi(x,y,\phi,\phi_x,\phi_{xx},\dots,\phi_y,\phi_{yy},\dots)
\]
of the Moutard equation~\eqref{eq:VsG}.
The resulting Voss nets are essentially those that arise along the
path followed by Guichard~\cite[Sections~V--XII]{CG-1890b}.
We start with some important cases of degeneracy,
 understood in the sense of Proposition~\ref{VsG->Voss}\,(ii),
that is, \cite[Theorem 2.25]{Izm}.\footnote{In
Guichard~\cite{CG-1890b}, degeneracy is expressed by saying that
``surfaces $S$ are spheres.''}

\begin{Proposition}[{Guichard~\cite[Section~XII]{CG-1890b}}]\label{prop:degen:tw}
The Voss surface associated to $\Phi = \phi_x$ is degenerate
if and only if the initial sine-Gordon solution satisfies the travelling
wave ansatz $\phi = \phi(x + m y)$, $m = {\rm const}$.
\end{Proposition}

\begin{proof} ``$\Rightarrow$'':
Assuming degeneracy, $\Phi = \phi_x$ satisfies $X = Y = 0$ by
Proposition~\ref{VsG->Voss}\,(ii),
which leads to
\[
\phi_{yy} = \frac{\phi_y}{\phi_x} \sin\phi, \qquad
\phi_{xy} = \sin \phi, \qquad
\phi_{xx} = \frac{\phi_x}{\phi_y} \sin\phi.
\]
The first equation is obtained by inserting $\Phi = \phi_x$ into
$Y$ given by formula~\eqref{eq:Phi:XY},
the second one is the sine-Gordon equation, while the third one
is the compatibility condition between the first two.
It follows that
$(\phi_y/\phi_x)_x = 0 = (\phi_y/\phi_x)_y$, that is,
$\phi_y/\phi_x = m$ = const.
Hence, $\phi = \phi(x + m y)$,
which is the travelling wave ansatz.

``$\Leftarrow$'':
Let $\phi_y = m \phi_x$.
Then $X = Y = 0$ by straightforward computation.
\end{proof}

By linearity, Proposition~\ref{prop:degen:tw} extends
to all translational symmetries $a \phi_x + b \phi_y$.
Thus, translational symmetries produce only degenerate Voss surfaces
from both the Dini helicoids and the Beltrami pseudospheres.
Still, one can obtain non-trivial results by using the scaling
symmetry.

\begin{Example}[the {\it right helicoid}]
\label{ex:rhelicoid}
Applying the scaling symmetry $-\frac12 x \phi_x + \frac12 y\phi_y$
to the unit-velocity kink
$\phi = 4 \arctan(\mathrm{e}^{x + y})$, we get
\[
\Phi = -\frac x2 \phi_x + \frac y 2 \phi_y
 = -\frac{x - y}{\cosh(x + y)}.\]
To obtain the associated Voss surface, we
insert $\Phi$ into formula~\eqref{eq:Voss:Phi} along with
$\ve r$ and $\ve n$ given by equations~\eqref{eq:PS:r} and~\eqref{eq:PS:n}, respectively.
The resulting Voss net is
\[
\ve q = \left[-\frac{\sin(x - y)}{\sinh(x + y)}, \frac{\cos(x - y)}{\sinh(x + y)},
x - y\right],
\]
that is, the common right helicoid in one of its Voss parameterisations,
see Figure~\ref{fig:helicoid}.
\begin{figure}[t]\centering
\includegraphics[scale=0.70,angle=90]{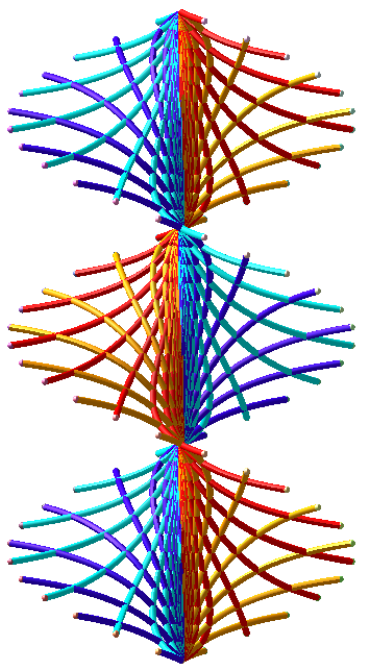}
\caption{The Voss net on the right helicoid from the pseudosphere and the scaling symmetry.}\label{fig:helicoid}
\end{figure}
\end{Example}

\begin{Example}[\emph{Kostin's helicoid}]
\label{ex:Dhelicoid}
Applying the scaling symmetry to a~kink of arbitrary speed and the
corresponding Dini helicoid~\eqref{eq:Dini:r}, we get
\begin{equation}
\label{eq:Dini2Voss}
\ve q = \left[-\frac{\sin v}{\sinh u},
 \frac{\cos v}{\sinh u},
 v + (u - \coth u) \sinh\gamma\right],
\end{equation}
where $u = \mathrm{e}^\gamma x + \mathrm{e}^{-\gamma} y$ and $v = x - y$;
it collapses to the right helicoid when $\gamma = 0$.
The curvatures are{\samepage
\[
K = -\frac{\tanh^4 u}{4 \cosh^2\gamma} ,
\qquad
H = -\frac12 \sinh u \tanh\gamma.
\]
For pictures, see Figure~\ref{fig:Dini2Voss}.}

The surface turns out to be one of the Minkowskian analogues of the Dini
helicoid, found by Kostin \cite[Theorem~3]{AVK-2020}.
The generating curve
coincides with a~(suitably reparameterised) Minkowskian analogue of the
tractrix \cite[equation~(3)]{AVK-2020}.
For every tangent line, the Lorentz distance between the
point of tangency and the intersection with the $\eta$-axis is $\sqrt{-1}$,
if the ``speed of light'' is $1/\sqrt{\sinh\gamma}$.
\begin{figure}[t]\centering
\includegraphics[scale=0.67,angle=90]{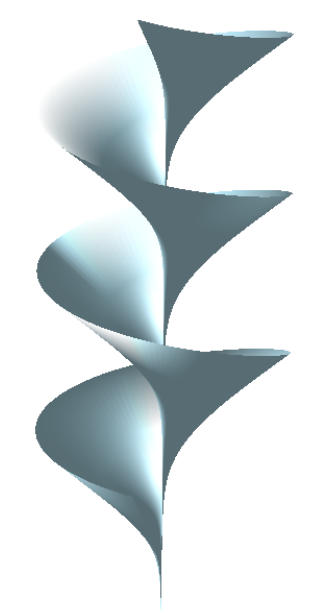} \qquad
\includegraphics[scale=0.49,angle=90]{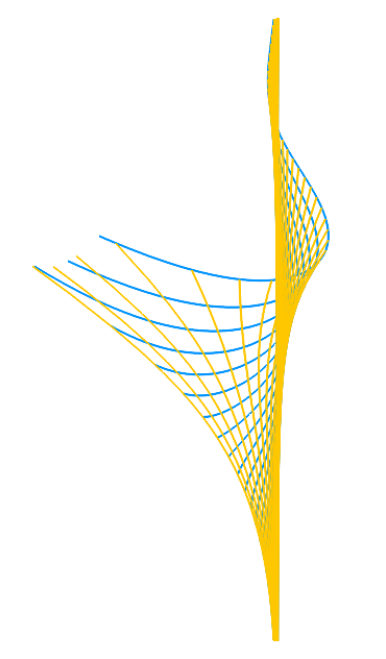}
\caption{Kostin's helicoid (left); a~part of its Voss net (right).}
\label{fig:Dini2Voss}
\end{figure}
\end{Example}

Turning back to the translational symmetries, they are still
applicable to other pseudospherical surfaces.
We provide only one simple example.

\begin{Example}[the {\it koru surface}]
\label{ex:Voss:koru}
Let us apply translational symmetries to the degenerate two-soliton
sine-Gordon solution
\begin{equation}
\label{eq:Kuen:phi}
\phi = 4 \arctan\left(\frac{x - y}{\cosh(x + y)}\right),
\end{equation}
see~\cite[Section~2.6.2]{P-M-2007} or~\cite[equation~(2.12)]{JC-FC-AF-1917}.
The corresponding pseudospherical surface is the Kuen surface
\[
\ve r = \frac{2 \cosh u}{v^2 + \cosh^2 u}
 [\cos v + v \sin v, \sin v - v \cos v, -{\sinh u}]
 + [0,0,u],
\]
where $u = x + y$, $v = x - y$.
The unit normal is
\[
\ve n = \frac{2 v}{v^2 + \cosh^2 u} [
-{\cos v}, -{\sin v}, \sinh u]
+ \frac{v^2 - \cosh^2 u}{v^2 + \cosh^2 u}
 [-{\sin v}, \cos v, 0].
\]

For $\Phi = \phi_x + \phi_y$, the associated Voss
surface~\eqref{eq:Voss:Phi} degenerates to a~point.
Therefore, we lose nothing if we restrict ourselves to the multiples
of~$\phi_x$.

The translational symmetry $\frac14 \phi_x$ maps the Kuen
surface to the Voss net
\[
\left[\frac{v \sin v - \sinh^2 u \cos v}
 {v \cosh^3 u - v^3 \cosh u},
- \frac{v \cos v + \sinh^2 u \sin v}
 {v \cosh^3 u - v^3 \cosh u},
\frac12 \frac{\sinh^3 u - \bigl(v^2 + 1\bigr) \sinh u}
 {v \cosh^3 u - v^3 \cosh u}\right].
\]
The real singular points where $\det\mathrm{I} = 0$ (ridges)
are determined by
$
2 v = \pm \sinh 2u$,
while the blow-up points where $\det\mathrm{I} = \infty$ (points at infinity)
are determined by one of
$v = 0$, $ v = \pm \cosh u$.
The curves intersect (see Figure~\ref{fig:koru:sing}),
where the ridges extend to infinity.

\begin{figure}[!ht]\centering
\unitlength=1mm
\begin{picture}(90,40)(0,-3)
\put(5,5){\includegraphics[scale=0.4,angle=90]{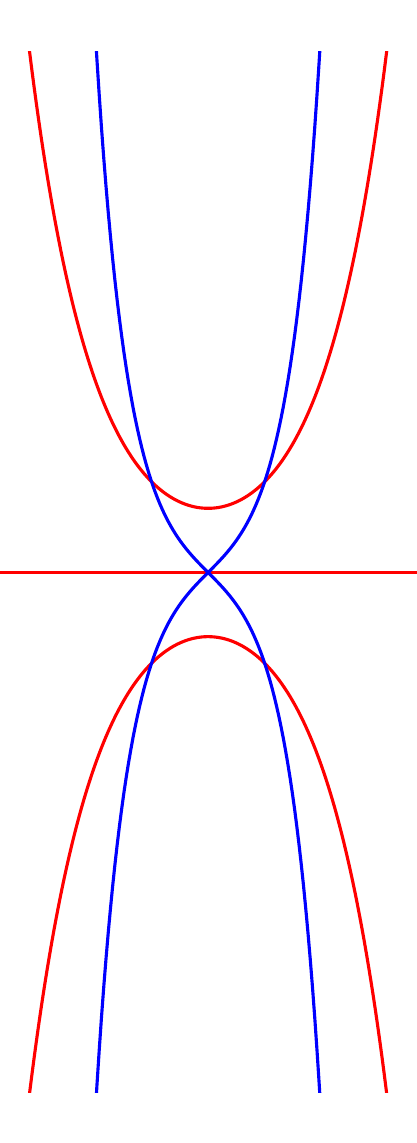}}
\put(0,0){\vector(1,0){85}}
\put(43.8,-1){\line(0,1){2}}
\put(43.8,-4){\hbox to 0pt{\small\hss$0$\hss}}
\put(80,-4){\hbox to 0pt{\small\hss$v$\hss}}
\put(0,0){\vector(0,1){35}}
\put(-1,19.1){\line(1,0){2}}
\put(-4,18){\hbox to 0pt{\small\hss$0$\hss}}
\put(-4,31){\hbox to 0pt{\small\hss$u$\hss}}
\end{picture}

\vspace{-1mm}

\caption{Singular and blow-up points of the koru surface
($\det\mathrm{I} = 0$ in blue,
$\det\mathrm{I} = \infty$ in red).}
\label{fig:koru:sing}
\end{figure}
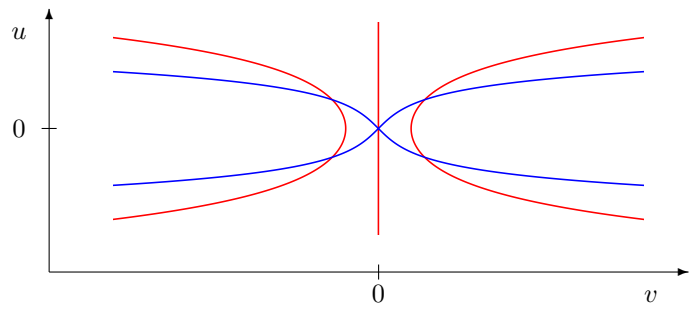

To explain the name given to the surface, in Figure~\ref{fig:koru}
we show part of the surface delimited by $|u| < 1$, $\cosh u < v$
next to the Hundertwasser koru flag of New Zealand~\cite{HKF}.

\begin{figure}[!ht]\centering
\includegraphics[scale=0.407]{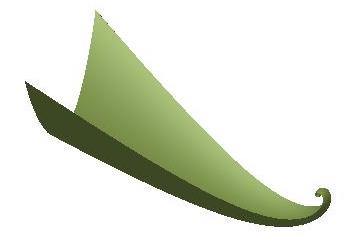} \qquad
\raise 9mm \hbox{\includegraphics[scale=0.11]{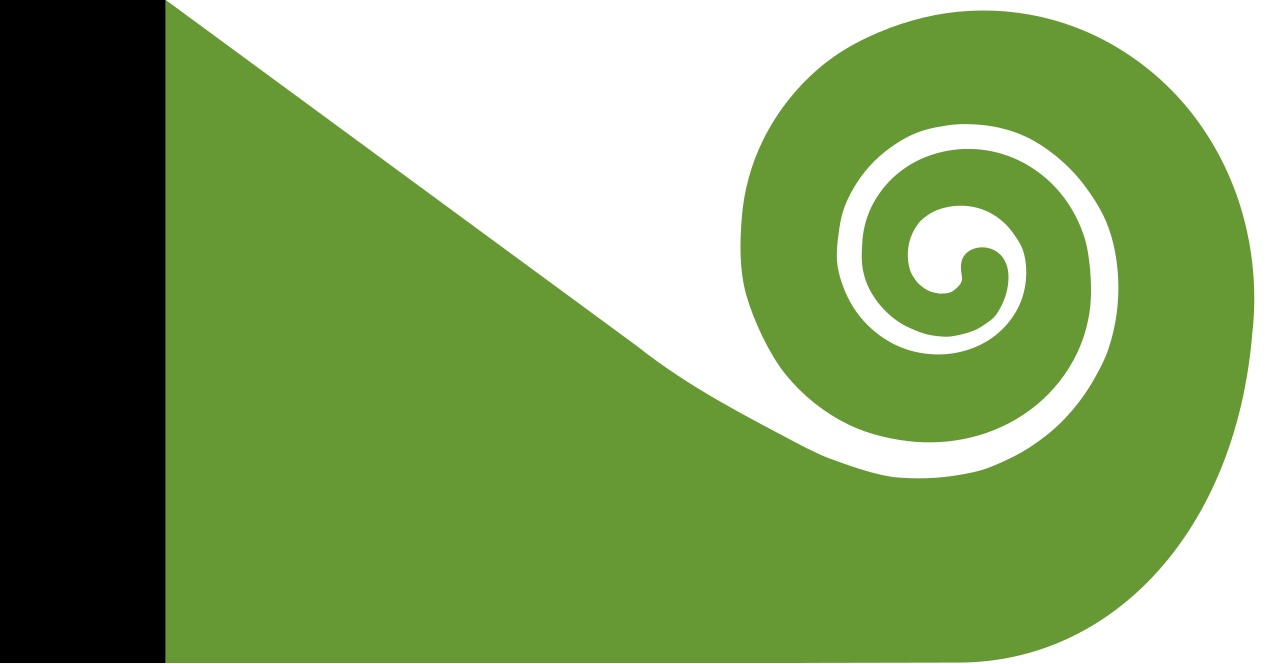}}
\vspace{-1mm}

\caption{Part of the koru surface next to the
Hundertwasser's koru flag.}\label{fig:koru}
\end{figure}

\begin{figure}[!ht]\centering
\unitlength = 1mm
\begin{picture}(150, 56)
\put(0, 20){\includegraphics[scale=0.3]{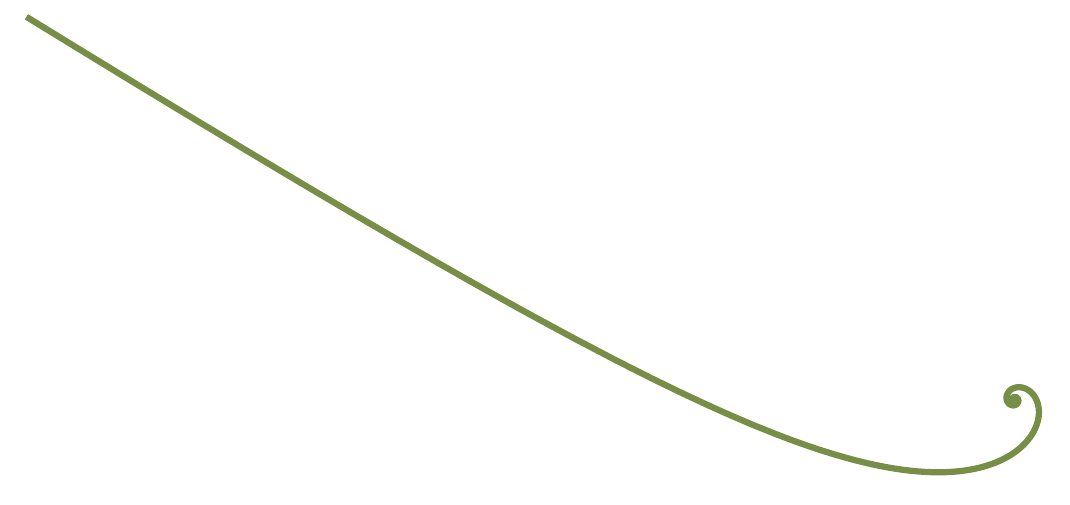}}
\put(66, 0){\includegraphics[scale=0.3]{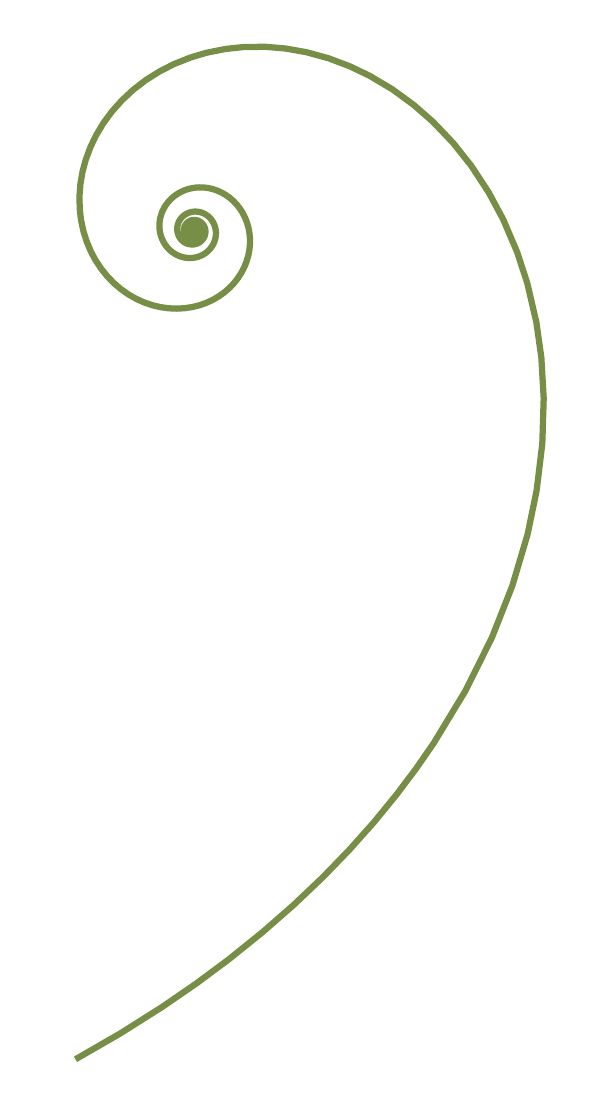}}
\put(98, 0){\includegraphics[scale=0.3]{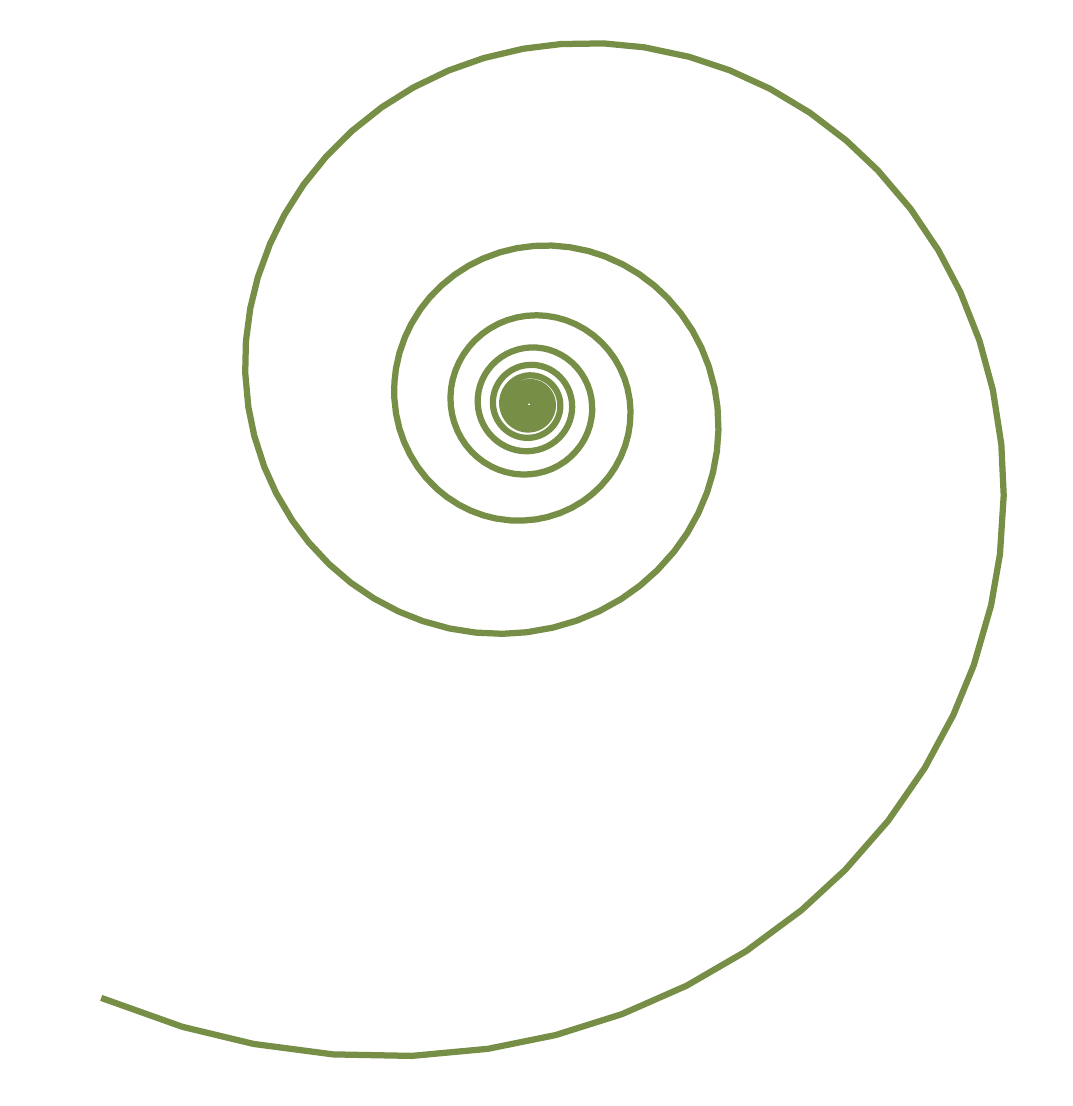}}
\color{Gray}
\thicklines
\put(50.6, 21.3){\framebox(4, 5.6){}}
\put(72.5, 40.5){\framebox(7.6, 8.3){}}
\end{picture}
\vspace{-1mm}

\caption{The infinitely curled spinal curve 
of the koru surface in three magnifications.}
\label{fig:fern}
\end{figure}
\end{Example}

\begin{table}[!ht]\centering
\begin{tabular}{r|rr}
 & translation & scaling
\\\hline
pseudosphere & degenerate \small (Proposition~\ref{prop:degen:tw})
 & right helicoid \small (Example~\ref{ex:rhelicoid})
\\
Dini helicoid & degenerate \small (Proposition~\ref{prop:degen:tw})
 & Kostin's helicoid \small (Example~\ref{ex:Dhelicoid})
\\
Kuen surface & koru surface \small (Example~\ref{ex:Voss:koru})
 & eared screw \small (Example~\ref{ex:Voss:archim})
\end{tabular}
\caption{Summary of examples.}
\label{tab:ex}
\end{table}

Resembling an unfurling fern frond,
the ``spinal curve'' $u = 0$ is the spiral $r = 1/\bigl(v^2 - 1\bigr)$,
that is, the inversion of the Galileo spiral $r = v^2 - 1$
\cite{VR-2000}.
Here $r$ and $v$ mean the radial coordinates (radius and angle).

\begin{Example}[the {\it eared screw}]
\label{ex:Voss:archim}
Applying the scaling symmetry to the same
Kuen surface as in Example~\ref{ex:Voss:koru}, we get the Voss net
\begin{gather*}
\biggl[A \sin v + B \cos v, B \sin v - A~\cos v,
\\
\qquad
\frac{v^4 \cosh u
 - v^2\bigl(\cosh^3 u + \cosh u + u \sinh u\bigr)
 + \bigl(1 - \sinh^2 u\bigr) (\cosh u - u \sinh u)}
 {v \cosh u \bigl(\cosh^2 u - v^2\bigr)}\biggr],
\\
A= \frac{2 u + \sinh 2 u}
 {\cosh u\bigl(\cosh^2 u - v^2\bigr)},
\qquad
B = \frac{u - u \cosh 2 u + \bigl(1 - v^2\bigr) \sinh 2 u}
 {v \cosh u \bigl(\cosh^2 u - v^2\bigr)}.
\end{gather*}
The singularity condition $\det\mathrm{I} = 0$ is equivalent to
\[
\bigl(v^2 - 1\bigr) \cosh 2 u + (u \pm v) \sinh 2 u + v^2 \mp 2 u v - 1 = 0,
\]
whereas $\det\mathrm{I} = \infty$ iff
$
v = \pm\cosh u$,
see Figure~\ref{fig:archim:sing}.
The intersections of blue and red lines lie at~${v = \pm \sqrt 2, \pm 1, 0}$.

\begin{figure}[!ht]\centering
\unitlength=1mm
\begin{picture}(90,40)(0,-3)
\put(8,1.6){\includegraphics[scale=0.5,angle=90]{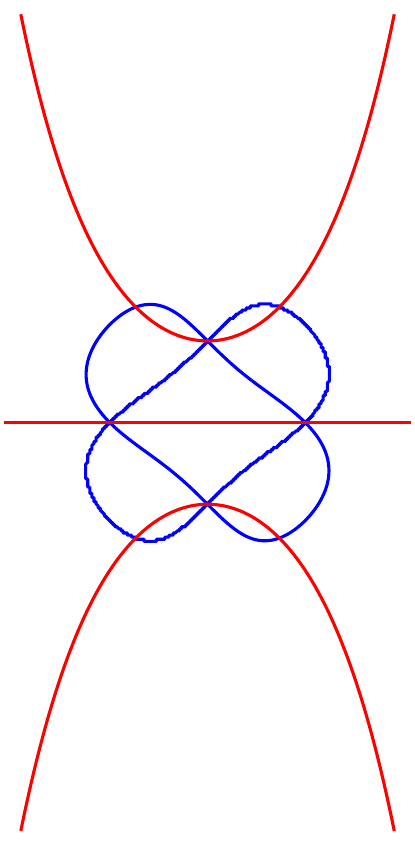}}
\put(0,0){\vector(1,0){85}}
\put(34,-1){\line(0,1){2}}
\put(36.85,-1){\line(0,1){2}}
\put(43.8,-1){\line(0,1){2}}
\put(50.75,-1){\line(0,1){2}}
\put(53.6,-1){\line(0,1){2}}
\put(43.8,-4.5){\hbox to 0pt{\small\hss$0$\hss}}
\put(51,-4.5){\hbox to 0pt{\small\hss$1$\hss}}
\put(31.5,-4.5){\hbox to 0pt{\small\hss$-\sqrt 2$\hss}}
\put(80,-4.5){\hbox to 0pt{\small\hss$v$\hss}}
\put(0,0){\vector(0,1){38}}
\put(-1,19.1){\line(1,0){2}}
\put(-4,18){\hbox to 0pt{\small\hss$0$\hss}}
\put(-4,31){\hbox to 0pt{\small\hss$u$\hss}}
\end{picture}
\caption{Singular points of the eared screw surface
($\det\mathrm{I} = 0$ in blue,
$\det\mathrm{I} = \infty$ in red).}
\label{fig:archim:sing}
\end{figure}

Although resembling the right helicoid, the surface is not minimal.
In the area $v > \cosh u$, the surface looks like a~screw with ears,
see Figure~\ref{fig:archim:ears}.
For the corresponding Voss net in a~different parameter domain, see
Figure~\ref{fig:archim:net}.

\begin{figure}[!ht]\centering
\unitlength=1mm
\begin{picture}(120,55)(-55,-3)
\put(-63,-2){\includegraphics[scale=0.45]{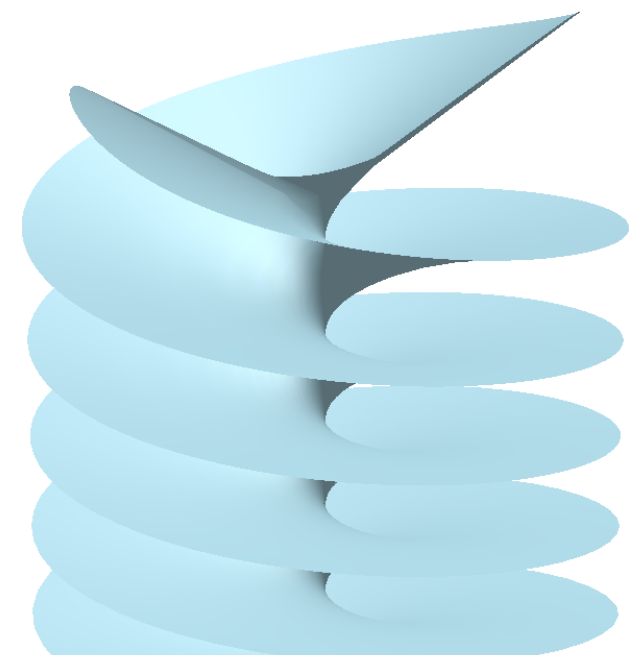}}
\put(28,46){\includegraphics[scale=0.4,angle=270]{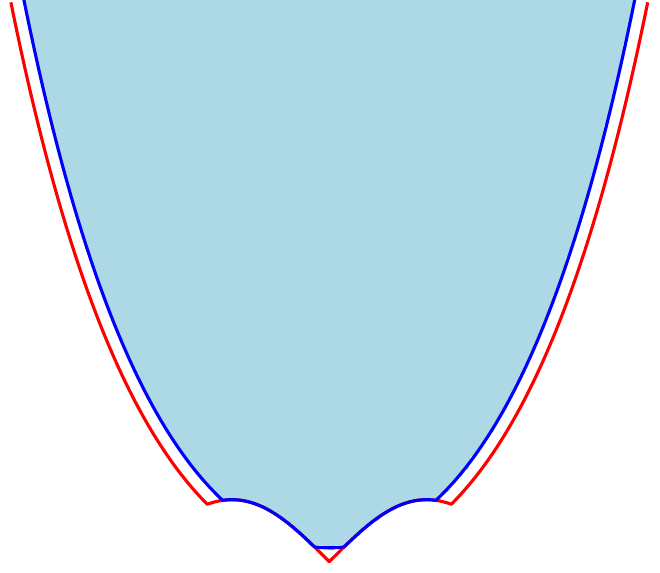}}
\put(0,0){\vector(1,0){70}}
\put(28.8,-1){\line(0,1){2}}
\put(28.8,-4.5){\hbox to 0pt{\small\hss$1$\hss}}
\put(64,-4.5){\hbox to 0pt{\small\hss$v$\hss}}
\put(0,0){\vector(0,1){46}}
\put(-1,23.7){\line(1,0){2}}
\put(-4,22.7){\hbox to 0pt{\small\hss$0$\hss}}
\put(-4,40){\hbox to 0pt{\small\hss$u$\hss}}
\end{picture}
\caption{The eared screw surface (left)
 and its parameter domain (right).}
\label{fig:archim:ears}
\end{figure}

\begin{figure}[!ht]\centering
\unitlength=1mm
\begin{picture}(120,59)(-62,-5)
\put(-58,-6){\includegraphics[scale=0.55]{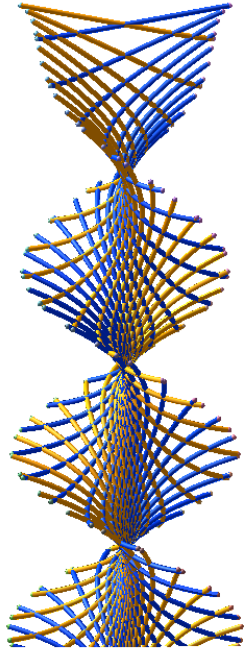}}
\put(0,0){\includegraphics[scale=0.34,angle=0]{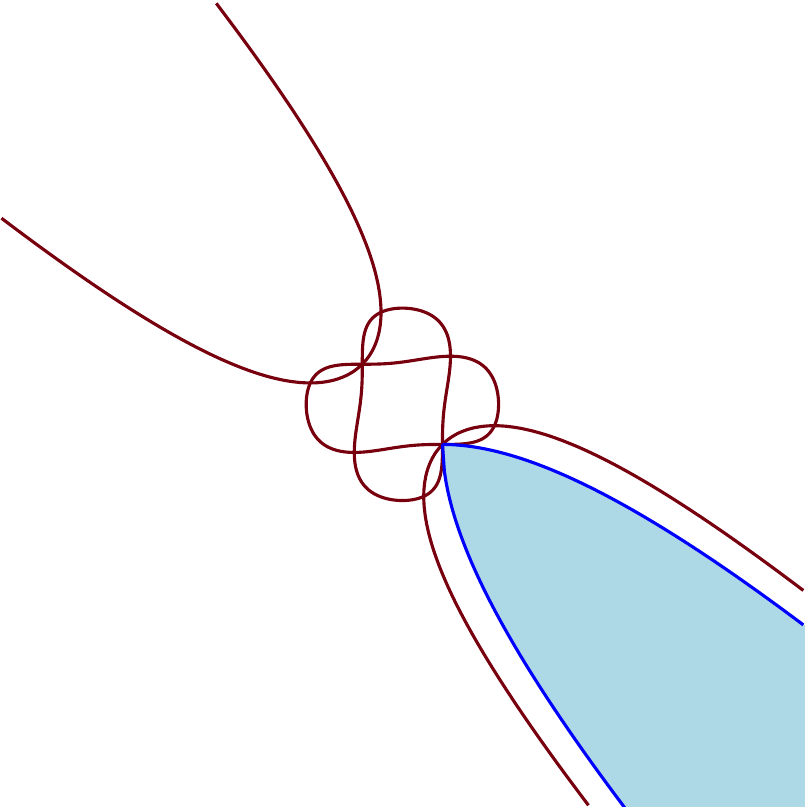}}
\put(0,0){\vector(1,0){51}}
\put(23.1,-1){\line(0,1){2}}
\put(23.1,-4.5){\hbox to 0pt{\small\hss$0$\hss}}
\put(46,-4.5){\hbox to 0pt{\small\hss$x$\hss}}
\put(0,0){\vector(0,1){47}}
\put(-1,23.2){\line(1,0){2}}
\put(-2,22){\hbox to 0pt{\small\hss$0$}}
\put(-2,42){\hbox to 0pt{\small\hss$y$}}
\end{picture}
\caption{The Voss net of the eared screw surface (left)
 and its parameter domain (right).}
\label{fig:archim:net}
\end{figure}

\end{Example}

\section{Recursion operators of the sine-Gordon equation}
\label{sect:sG:RO}

The prevailing viewpoint is that recursion operators are linear
operators, mapping symmetries to symmetries.
They are available for many equations and have a~rich and deep theory (see, e.g.,
\cite{ISK-PHMK-2000,JK-AV-RV-2017,S-W-2001b,VEZ-BGK-1984}
and references therein)
developed primarily in the integro-differential formalism,
relying on formal inverses of total derivatives.

The integro-differential recursion operator due to
Olver~\cite[equation~(21)]{PJO-1977}
acts on a~sine-Gordon symmetry $\Phi$ by the formula
\begin{equation}
\label{eq:sG:RO:Olver}
\RO\Phi = D_x^2 \Phi + \phi_x^2 \Phi - \phi_x D_x^{-1} (\phi_{xx} \Phi) = D_x^2 \Phi + \phi_x D_x^{-1} (\phi_x D_x \Phi).
\end{equation}
The latter version results from integration by parts.

As stressed by many authors
\cite{GAG-1994,PJO-1977,S-W-2001a,JMV-1991},
inverses of total derivatives are poorly defined.
Quite often (see~\cite{L-H-L-2021,S-W-2001a}),
the formal inverse $D_x^{-1}$ is introduced
by the contradictory\footnote{Indeed, $D_x^{-1} D_x f(y) = f(y)$
requires that $D_x 0 = f(y)$ for every function $f(y)$.\label{foo:Dinv}}
rule~${D_x \circ D_x^{-1} = D_x^{-1} \circ D_x = \mathrm{Id}}$.
Moreover, although a~formal integro-differential inverse recursion operator
is available for~$\RO$, see~\cite[Section~3]{Ai-1983}, \cite{K-K-2002,L-C-1993}, it is insufficient for our purposes,
mainly because we shall need to invert also $\RO + \mathrm{Id}$,
which is something very different.

This is why we adhere to Guthrie's~\cite{GAG-1994} formalism,
treating recursion operators as B\"acklund auto-transformations
for the linearised equation~\cite{H-Kh-P-2016,MM-alro,MM-rzcr,CJP-1991}.
If, following this path, we identify the linearised sine-Gordon equation with
the Moutard equation, we obtain Proposition~\ref{prop:sG:RO:Guthrie} below,
which, however, merely reproduces, up to notational variance, a~result published by Guichard in 1890.\looseness-1

\begin{Proposition}[{\cite[equations~(3) and (4)]{CG-1890a},
\cite[Section~V]{CG-1890a}; \cite[Example~4.4]{MM-alro}}]
\label{prop:sG:RO:Guthrie}
 Given a~solution~$\phi$ of the sine-Gordon equation,
let~$\Phi$ be a~solution of the Moutard equation.
Then the system
\begin{equation}
\label{eq:VsG:tilde}
Q_x = \Phi_x \phi_x,
\qquad
Q_y = \Phi \sin \phi
\end{equation}
is compatible and
$
\Psi = \Phi_{xx} + \phi_x Q
$
is another solution of the
Moutard equation.
\end{Proposition}

\begin{proof}
Equations~\eqref{eq:VsG:tilde} pass the cross-differentiation test
\[
Q_{xy} = (\Phi_x \phi_x)_y
 = \Phi_{xy} \phi_x + \Phi_x \phi_{xy}
 = \Phi \phi_x \cos\phi + \Phi_x \sin\phi
 = (\Phi \sin \phi)_x
 = Q_{yx},
\]
while $\Psi$ satisfies
\begin{gather*}
\Psi_{xy} = \Phi_{xxxy} + (\phi_x Q)_{xy}
 = (\Phi \cos\phi)_{xx}
 + \phi_{xxy} Q + \phi_{xx} Q_y + \phi_{xy} Q_x + \phi_x Q_{xy}
\\
\phantom{\Psi_{xy}}{}
 = \Phi_{xx} \cos\phi + \phi_x Q \cos\phi = \Psi \cos\phi,
\end{gather*}
which finishes the proof.
\end{proof}

To see that Proposition~\ref{prop:sG:RO:Guthrie} actually delivers
the Guthrie form of Olver's operator~\eqref{eq:sG:RO:Olver},
it is sufficient to observe that $Q$ can be chosen for
$D_x^{-1} (\phi_x D_x \Phi)$.
Thus, instead of~\eqref{eq:sG:RO:Olver} we can write
\begin{gather}
\label{eq:sG:RO:Guthrie}
\RO\Phi = \Phi_{xx} + \phi_x Q,
\qquad
Q_x = \Phi_x \phi_x,
\qquad
Q_y = \Phi \sin \phi.
\end{gather}

Note that the potential $Q$ coincides
with the linearisation of the potential $q_1$ associated with one of the
conservation laws of the
sine-Gordon equation, see formula~\eqref{eq:q} below.
This origin of $Q$ was to be expected, cf.\ \hbox{\cite[Section~4]{MM-alro}}.

\begin{Remark}
To emphasise the multivalued character of the operator $\RO$,
Proposition~\ref{prop:sG:RO:Guthrie} can be illustrated by the diagram
\[
\qquad\begin{array}{ccccc}
&&\fpa{Q_x &= \Phi_x \phi_x,
\\
Q_y &= \Phi \sin \phi
\\
\Phi_{xy} &= \Phi \cos\phi
\\
\phi_{xy} &= \sin\phi}{2.5cm}
\\
& \raise1ex\llap{\scriptsize} \swarrow && \searrow
\raise1.4ex\rlap{\hglue-1.2ex\scriptsize $\Psi = \Phi_{xx} + \phi_x Q$}
\\
\fpa{\Phi_{xy} &= \Phi \cos\phi
\\
\phi_{xy} &= \sin\phi}{2.5cm}
&&&&\fpa{\Psi_{xy} &= \Psi \cos\phi
\\
\phi_{xy} &= \sin\phi}{2.5cm}
\end{array}
\]
The arrow on the left forgets $Q$.
Every solution $\Phi $ of the Moutard equation on the left
has a~one-parameter family of preimages $(\Phi,Q)$,
which is mapped to a~one-parameter family of solutions~$\RO\Phi$
of the Moutard equation on the right.
The parameter is the integration constant for $Q$.
\end{Remark}

\begin{Remark}
\label{rem:multi}
In the Guthrie formalism,
the fact that recursion operators are multi-valued cannot be ignored.
Yet one commonly writes $\Psi = \RO\Phi$ with due reservations
and caution if $\Psi$ is among
the results for a~specific choice of the integration constant.
When we want to stress the multi-valuedness, we write
$\Psi \equiv \RO\Phi$ or,
more accurately, $\Psi \equiv \RO\Phi \operatorname{mod} \RO 0$,
where $\RO 0 = \phi_x \RR$
is the $\RR$-ideal consisting of all results obtained from
zero for all possible values of the integration constant.
\end{Remark}

Since Guthrie's recursion operators are essentially relations,
the natural way to introduce their inverses is by the equivalence
\begin{equation}
\label{eq:RO:inv}
\Phi \equiv \RO^{-1}\Psi
\quad\stackrel{\rm def}{\Leftrightarrow}\quad
\Psi \equiv \RO\Phi.
\end{equation}
Alternatively speaking, to obtain $\RO^{-1}$ in the Guthrie form,
all that is needed is to express $\Phi$ in terms of~$\Psi$,
see \cite[Section~4]{GAG-1994}; for examples see also \cite[Section~6]{MM-alro},
\cite[Section~6.1]{MM-rzcr}
and Proposition~\ref{prop:sG:iRO:Guthrie} below.

When exchanging $x$ and $y$ in Proposition~\ref{prop:sG:RO:Guthrie},
we obviously get another recursion operator, which we shall denote by $\RO'$,
namely
\begin{equation}
\label{eq:RO'}
\RO'\Psi = \Psi_{yy} + \phi_y Q',
\qquad
Q'_y = \Psi_y \phi_y,
\qquad
Q'_{x} = \Psi \sin\phi.
\end{equation}
As proved already by Guichard in 1890,
$\RO'$ is inverse to $\RO$ up to terms induced by
the inherent ambiguity; this original version can be expressed by formulas
\eqref{eq:Gui:rel} below.
The integro-differential version \cite[Proposition~2]{L-H-L-2021} of this result
was first presented only in 2021,
without taking into account the ambiguities.

The ambiguity-resistant version looks as follows.

\begin{Proposition}
\label{prop:sG:iRO:Guthrie}
The inverse recursion operator~$\RO^{-1}$ defined by
\hbox{\eqref{eq:RO:inv}} coincides with the operator
$\RO'$ defined by~\eqref{eq:RO'}.
\end{Proposition}

\begin{proof}
According to Proposition~\ref{prop:sG:RO:Guthrie},
the equality $\Psi = \RO\Phi$ can be expressed by the compatible system
\begin{gather}
\phi_{xy} = \sin\phi,
\qquad
\Phi_{xy} = \Phi \cos\phi,
\qquad
\Psi = \Phi_{xx} + \phi_x Q,
\nonumber\\
Q_x = \Phi_x \phi_x,\qquad
Q_y = \Phi \sin \phi.\label{eq:RO:sys}
\end{gather}
To express $\Phi$ in terms of~$\Psi$,
one can apply the standard elimination-completion algorithm,
see, e.g., \cite{R-L-W-2001} and references therein,
under a~suitable elimination ordering of derivatives.
The result will be an equivalent system, possessing the same solutions.
To start with, we express $Q$ from the third equation~\eqref{eq:RO:sys},
obtaining $Q = (\Psi - \Phi_{xx})/\phi_x$.
Inserting into the last two equations, we obtain
\[
\Phi_{xx} = \frac{\phi_x}{\sin\phi} (\Phi_x \cos\phi - \Psi_y) + \Psi,
\qquad
\Phi_{xxx} = \frac{\phi_{xx}}{\sin\phi} (\Phi_x \cos\phi - \Psi_y)
 - \Phi_x \phi_x^2 + \Psi_x.
\]
Evaluating the compatibility conditions
$\Phi_{xxx} = D_x \Phi_{xx}$ and
$D_x \Phi_{xy} = D_y \Phi_{xx}$, we get
\[
\Psi_{xy} = \Psi \cos\phi
 \qquad \text{and} \qquad \Phi_x = \frac{\sin\phi}{\phi_y} (\Phi - \Psi_{yy}) + \Psi_y \cos\phi,
\]
respectively.
The third equation~\eqref{eq:RO:sys} is thereby satisfied,
while the fourth yields
\[
\Phi_y = \Psi_{yyy} + \frac{\phi_{yy}}{\phi_y} (\Phi - \Psi_{yy})
 + \Psi_y \phi_y^2.
\]
Thus, we are left with the system
\begin{gather*}
\phi_{xy} = \sin\phi,
\qquad
\Psi_{xy} = \Psi \cos\phi,
\qquad
\Phi_x = \frac{\sin\phi}{\phi_y} (\Phi - \Psi_{yy}) + \Psi_y \cos\phi,
\\
\Phi_y = \Psi_{yyy} + \frac{\phi_{yy}}{\phi_y} (\Phi - \Psi_{yy})
 + \Psi_y \phi_y^2,
\end{gather*}
which is compatible now and equivalent to~\eqref{eq:RO:sys}.
Substituting $\Phi = \Psi_{yy} + \phi_y Q'$, we obtain another equivalent
system
\begin{gather*}
\phi_{xy} = \sin\phi,
\qquad
\Psi_{xy} = \Psi \cos\phi,
\qquad
\Phi = \Psi_{yy} + \phi_y Q',
\qquad
Q'_y = \Psi_y \phi_y,
\qquad
Q'_{x} = \Psi \sin\phi.
\end{gather*}
These are exactly formulas~\eqref{eq:RO:sys}
with $x$, $y$ and $\Phi$, $\Psi$ exchanged.
Therefore, $\Phi = \RO'\Psi$, proving that $\RO' = \RO^{-1}$ in the sense
of the equivalence~\eqref{eq:RO:inv}.
\end{proof}

Turning back to the original version published by
Guichard~\cite{CG-1890b,CG-1890a},
it can be stated as
\begin{equation} \label{eq:Gui:rel}
\RO' \RO \Phi \equiv \Phi \operatorname{mod} \RO'0, \qquad
\RO \RO' \Psi \equiv \Psi \operatorname{mod} \RO0
\end{equation}
in the notation of Remark~\ref{rem:multi}.
In Proposition~\ref{prop:sG:RO:Guthrie:inv} below,
we give this result yet another form,
which will be useful later.
First, we introduce trivial recursion operators as a~tool to tackle the
ambiguities.

\begin{Definition}\label{def:RO:triv}
By a~{\it trivial}~recursion operator we mean the operator
of the form
\[
\Z\Phi = W^{(1)} \Sigma^{(1)} + \cdots + W^{(n)} \Sigma^{(n)},
\]
where $\Sigma^{(i)}$ are fixed symmetries and $W^{(i)}$ satisfy
\smash{$
D_x W^{(i)} = D_y W^{(i)} = 0$}
(that is, $W^i = $ const).
\end{Definition}

Thus, $\Z$ formally satisfies Guthrie's definition
($\Z$ is formally a~B\"acklund transformation),
but $\Z\Phi$ does not depend of $\Phi$.

\begin{Proposition}[cf.\ {Guichard~\cite{CG-1890b,CG-1890a}}]
\label{prop:sG:RO:Guthrie:inv}
For~$\RO'$ defined by~\eqref{eq:RO'},
both $\RO' \circ \RO - \mathrm{Id}$ and~$\RO \circ \RO' - \mathrm{Id}$ are trivial.
\end{Proposition}

\begin{proof}
We have
$
\RO' \RO \Phi = (\RO \Phi)_{yy} + \phi_y Q'$,
where $Q'$ satisfies the last two equations~\eqref{eq:RO'},
that is,
$
Q^{\prime}_x = (\Phi_{xx} + \phi_x Q) \sin \phi$,
$
Q^{\prime}_y
 = (\Phi_x + Q \sin\phi) \phi_y$.
Denoting $W = Q' + Q \cos \phi - \Phi_x \sin \phi$, we have
$
D_x W = 0$, $D_y W = 0$.
However,
\begin{gather*}
\RO' \RO \Phi = (\Phi_{xx} + \phi_x Q)_{yy} + \phi_y Q'
\\
\phantom{\RO' \RO \Phi}{}
 = \Phi_{xxyy} + \phi_{xyy} Q + 2 \phi_{xy} Q_y + \phi_x Q_{yy}
 + (\Phi_x \sin \phi - Q \cos \phi + W) \phi_y
 = \Phi + W\phi_y
\end{gather*}
by~\eqref{eq:sG}, \eqref{eq:VsG} and~\eqref{eq:VsG:tilde}.
Thus, $\RO' \circ \RO - \mathrm{Id}$ is trivial.
At the same time, we obtained the first Guichard's relation~\eqref{eq:Gui:rel}.
Triviality of
$\RO \circ \RO' - \mathrm{Id}$ as well as the second Guichard's relation
follow by the $x \leftrightarrow y$ symmetry.
\end{proof}

 Following Olver~\cite{PJO-1977}, we construct an
infinite hierarchy of symmetries $\Phi_n$
(and thus an infinite number of solutions to the Moutard equation)
by the repeated application of the recursion operator.
In Guthrie's formalism, the procedure is
$
\Phi_{n + 1} = D_{xx} \Phi_n + \phi_x Q_n$,
where $Q_n$ satisfy~${
Q_{n,x} = \Phi_{n,x} \phi_x}$,
$
Q_{n,y} = \Phi_n \sin \phi$.
Depending on what we start with, we obtain two distinct hierarchies.
\begin{itemize}\itemsep=0pt
\item[(1)]
{\it The pmKdV hierarchy {\rm\cite[{\it Example}~6]{PJO-1977}}.}
Starting with $\Phi_1 = \phi_x$, we get
\begin{gather*}
Q_1 = \frac12 \phi_x^2,
\\
Q_2 =\phi_x \phi_{xxx} -\frac12 \phi_{xx}^2 + \frac38 \phi_x^4,
\\
\qquad\vdots
\end{gather*}
whence
\begin{gather}
\Phi_1 = \phi_x,
\nonumber\\
\Phi_2 = \phi_{xxx} + \frac12 \phi_x^3,
\nonumber\\
 \Phi_3
 = \phi_{xxxxx} + \frac52 \phi_x^2 \phi_{xxx}
 + \frac52 \phi_x \phi_{xx}^2 + \frac38 \phi_x^5,
\label{eq:hierarchy:pmKdV}
\\
\qquad\vdots
\nonumber
\end{gather}

\item[(2)] {\it The Khor'kova hierarchy~{\rm\cite[{\it Example}~2]{NGK-1988}}. }
Starting with the scaling symmetry
$\bar\Phi_1 = x \phi_x -\allowbreak y \phi_y$, we get
\begin{gather*}
\bar Q_1 = x Q_1 + q_1 + y \cos \phi,
\\
\bar Q_2 = x Q_2
 + q_1 Q_1
 - 3 q_2
 + 2 \phi_x \phi_{xx},
\\
\qquad\vdots
\end{gather*}
where the nonlocal variables $q_i$ are to be determined from{\samepage
\begin{alignat}{3}
& q_{1,x} = \frac{1}{2} \phi_x^2,
 \qquad &&
q_{1,y} = -{\cos \phi}, &\nonumber\\
& q_{2,x} = \frac{1}{2} \phi_{xx}^2 - \frac18 \phi_x^4,
 \qquad &&
q_{2,y} = \frac12 \phi_x^2 \cos \phi,&\label{eq:q}
\\
& \qquad\vdots && \qquad\vdots&\nonumber
\end{alignat}}%

\pagebreak

\noindent
as potentials of the standard conservation laws
\cite[equation~(1.5)]{D-B-1977}.
Hence the nonlocal hierarchy{\samepage
\begin{gather}
\bar\Phi_1 = x \phi_x - y \phi_y,
\nonumber\\
\bar\Phi_2
 = x \Phi_2 + q_1 \Phi_1 + 2 \phi_{xx},
\nonumber\\
\bar\Phi_3 =
 x \Phi_3
 + q_1 \Phi_2
 - 3 q_2 \Phi_1
 + 4 \phi_{xxxx}
 + 7 \phi_x^2 \phi_{xx},
\label{eq:hierarchy:NGK}
\\
\qquad\vdots
\nonumber
\end{gather}}%
The expressions $\bar\Phi_i$ represent {\it shadows} of nonlocal
symmetries in the sense of~\cite{B-V-V,ISK-AMV-1989,JK-AV-RV-2017},
meaning that they satisfy equation~\eqref{eq:VsG:D},
but depend also on the nonlocal variables $q_i$.
All of them can be extended to full symmetries in a~possibly larger
covering, see op.\ cit.
\end{itemize}

\begin{Remark} \label{rem:negative}
The local symmetries $\Phi_n$ and the nonlocal symmetry shadows
$\bar\Phi_n$ are referred to as members of the positive pmKdV hierarchy
and the positive Khor'kova hierarchy, respectively.
Using the inverse recursion operator $\RO'$ and starting with
$\Phi_{-1} = \phi_y$ and $\bar\Phi_{-1} = -\Phi_1$,
we get members of the corresponding negative hierarchies,
which will be denoted $\Phi_{-n}$ and \smash{$\bar\Phi_{-n}$}.
The explicit formulas are obtained by exchanging
$x \leftrightarrow y$ and $n \leftrightarrow -n$ in~\eqref{eq:hierarchy:pmKdV},
\eqref{eq:q}, and~\eqref{eq:hierarchy:NGK};
cf.\ Proposition~\ref{prop:sG:RO:Guthrie:inv}.
\end{Remark}

\begin{Remark} \label{rem:bt:i}
In the pmKdV and Khor'kova hierarchy above,
$Q_n$, $q_n$, $\bar Q_n$ have been chosen as homogeneous\footnote{Here, an expression $F(x,y,\phi,\phi_x,\phi_y,\dots)$ is said
to be homogeneous of weight $N$ if $F$ transforms to $\lambda^N F$
under $x \mapsto x/\lambda$,
$y\mapsto \lambda y$,
$\phi \mapsto \phi$.
See, e.g., \cite[Section~3.2]{JK-AV-RV-2017}.\label{fn:homo}}
of weights $2 n$, $2 n - 1$, $2 n - 1$, respectively, with regard to
weights $[x] = -1$, $[y] = 1$, $[\phi] = 0$,
which ensures that $\Phi_n$ and \smash{$\bar\Phi_n$} are homogeneous of
weight $2 n - 1$ and $2 n$, respectively.
This may be viewed as setting all integration constants to zero,
which is a~common practice, but not always a~well-defined procedure.
\end{Remark}

\section{Guichard sequences of Voss surfaces}
\label{sect:pos}

In the previous section, we recalled the infinite pmKdV hierarchy
of local symmetries $\Phi_n$ and the infinite Khor'kova hierarchy
of nonlocal symmetry shadows \smash{$\bar\Phi_n$},
for~$n$ both positive or negative.
Using them in the way described in Proposition~\ref{VsG->Voss},
we obtain families of Voss surfaces, which will be referred to as
{\it Guichard sequences}, although Guichard himself considered only
$\Phi_n$ and $\Phi_{-n}$.
We say that such a~sequence stops, if the resulting Voss
surfaces are degenerate or differ from the previous ones by a~shift
or correspond to the sum of previously obtained Voss surfaces.
This is so, for instance, for all travelling wave initial solutions
and the whole pmKdV hierarchy~$\Phi_n$,
as observed already by Guichard.

\begin{Proposition}[{\cite[p.~264]{CG-1890b}}]
\label{prop:tw:f_x}
All symmetries $\Phi_n$, $n \ne 0$, map every travelling wave
solution of the sine-Gordon equation to a~degenerate Voss surface.
\end{Proposition}

\begin{proof}
The travelling wave solutions satisfy
$\phi_y = m \phi_x$, where $m$ = const $\ne 0$.
Then
\[
\phi_{xx} = \frac{\sin\phi}{m}
\]
in consequence of the sine-Gordon equation.
Since $m \phi_x^2 + 2 \cos\phi$ is a~first integral,
the travelling wave solutions satisfy $m \phi_x^2 + 2 \cos\phi = c =$
const.

Starting with $\Phi_1 = \phi_x$, we get
\[
\Phi_2 = \phi_{xxx} + \frac12 \phi_x^3
 = \frac{\cos\phi}{m} \phi_x + \frac{c - 2\cos\phi}{2m} \phi_x
 = \frac{c}{2}\phi_x
 = \frac{c}{2}\Phi_1.
\]
Therefore, $\Phi_n = (c/2)^n \Phi_1$ by induction.
However, $\Phi_1$ yields a~degenerate Voss surface by
Proposition~\ref{prop:degen:tw}.
This proves the statement for all $n > 1$.
Analogously for $n < -1$.
\end{proof}

\begin{Proposition} \label{prop:rhelicoid:R}
All symmetries $\bar\Phi_n$, $n \ne 0$, map the Beltrami
pseudosphere to one and the same right helicoid up to a~shift.
\end{Proposition}

\begin{proof}
For $n = 1$, this is Example~\ref{ex:rhelicoid}.
For $n = 2$, we have
$\bar\Phi_{2} - \bar\Phi_1 = 2/{\cosh(x + y)}$,
for which the Voss surface degenerates to the point $[0,0,-2]$.
For $n > 2$, we have $\bar\Phi_n = 2\bar\Phi_{n-1} - \bar\Phi_{n-2}$,
which can be directly verified for $n = 3$
and extended to $n > 3$ by applying $\RO$.
Thus, surfaces corresponding to various
$\bar\Phi_n$, $n > 1$, differ by a~shift along the axis.
Analogously for $n < -1$. \hglue -1pc \null
\end{proof}

\begin{Proposition} \label{prop:Dhelicoid:R}
All symmetries $\bar\Phi_n$, $n \ne 0$, map the Dini
helicoid~\eqref{eq:Dini:r} to one and the same Kostin helicoid
\eqref{eq:Dini2Voss} up to a~shift.
\end{Proposition}

\begin{proof}
Analogously to the proof of Proposition~\ref{prop:rhelicoid:R}.
For $n = 1$, this is Example~\ref{ex:Dhelicoid}.
For $n = 2$, we have
$\bar\Phi_{2} - \mathrm{e}^{2\gamma} \bar\Phi_1
 = 2 \mathrm{e}^{2 \gamma}/{\cosh(\mathrm{e}^\gamma x +\mathrm{e}^{-\gamma} y)}$,
for which the Voss surface degenerates to the point
$\big[0,0,-\mathrm{e}^{3 \gamma} - \mathrm{e}^{\gamma}\big]$.
For $n > 2$, we have
$\bar\Phi_n = 2 \mathrm{e}^{2 \gamma} \bar\Phi_{n-1}
 - \mathrm{e}^{4 \gamma}\bar\Phi_{n-2}$.
Analogously for $n < -1$.
\end{proof}

\section{On Guichard's problem}
\label{sect:GP}

At the end of his paper, Guichard~\cite{CG-1890b} asked under what
conditions the process of generating Voss surfaces stops after the
$n$-th step.
To provide an answer, we prove a~lemma, which links the problem with the
symmetry invariance.
In what follows, by $\RO 0$ we mean the image of the zero symmetry.
According to formula~\eqref{eq:sG:RO:Guthrie},
$\RO 0$ coincides with the ideal $\phi_x \RR$.

\begin{Lemma} \label{lemma1}
The Voss surface corresponding to a~solution $\Phi$ of the
Moutard equation by formula~\eqref{eq:Voss:Phi} is degenerate if and only if
\begin{equation}
\label{eq:degen:RO:1}
\RO\Phi + \Phi \in \RO 0.
\end{equation}
\end{Lemma}

\begin{proof}
By Proposition~\ref{VsG->Voss}\,(ii),
the surface 
is degenerate if and only if $X = 0$ and $Y = 0$.

The ``$\Rightarrow$'' part.
The equalities $X = 0$ and $Y = 0$, together with
the Moutard equation, yield the system
\begin{gather}
\Phi_{xx} + \Phi
 - \left(\frac{\Phi_x}{\tan\phi} + \frac{\Phi_y}{\sin\phi}\right) \phi_x = 0,
\qquad
\Phi_{xy} - \Phi \cos\phi = 0,
\nonumber
\\
\Phi_{yy} + \Phi
 - \left(\frac{\Phi_y}{\tan\phi} + \frac{\Phi_x}{\sin \phi}\right) \phi_y = 0.
\label{eq:degen:iR0}
\end{gather}
Assuming that~\eqref{eq:degen:iR0} holds, we easily find that
\[
Q = -\frac{\Phi_x}{\tan\phi}-\frac{\Phi_y}{\sin\phi} + \text{const},
\]
whence
\[
\RO\Phi = \Phi_{xx} + \phi_x Q = -\Phi + \text{const}\, \phi_x
\]
by formula~\eqref{eq:sG:RO:Guthrie}.

The ``$\Leftarrow$'' part.
If
$\Phi_{xx} + \phi_x Q = \RO\Phi = -\Phi + \text{const}\, \phi_x$, then
\[
Q = -\frac{\Phi_{xx} + \Phi}{\phi_x} + \text{const}.
\]
Equating $Q_y \!=\! \Phi \sin \phi$ according to
equation~\eqref{eq:VsG:tilde},
we get the first line of system~\eqref{eq:degen:iR0} by~straightforward
computation,
which means that $X = 0$.
Then also $Y = 0$ by
Proposition~\ref{VsG->Voss}\,(iii).
\end{proof}

\begin{Remark}
The condition~\eqref{eq:degen:RO:1} is equivalent to
\begin{equation}
\label{eq:degen:Bi}
\mathrm{III}^{ij} \Phi_{,i }\Phi_{,j} + \Phi^2 = \mathrm{const,}
\end{equation}
where
\[
\mathrm{III}^{ij} \Phi_{,i }\Phi_{,j}
 = \frac{\Phi_x^2 + 2 \Phi_x \Phi_y \cos\phi + \Phi_y^2}{\sin^2 \phi}
\]
is the Beltrami invariant of the Gaussian image of~$\Phi$.
This follows from the fact that the expression on the left-hand side
of equation~\eqref{eq:degen:Bi} is a~first integral of
system~\eqref{eq:degen:iR0}.
\end{Remark}

Turning back to the Guichard problem, we reword it as follows.
Given a~sine-Gordon solution~$\phi$, consider
members $\Phi_n$ of the pmKdV hierarchy, described by
formula~\eqref{eq:hierarchy:pmKdV},
members $\bar\Phi_n$ of the Khor'kova hierarchy, described by
formula~\eqref{eq:hierarchy:NGK},
and their negative counterparts~$\Phi_{-n}$ and~\smash{$\bar\Phi_{-n}$},
which correspond to~$\Phi_n$ and~\smash{$\bar\Phi_n$} by the discrete
symmetry $x \leftrightarrow y$ according to Remark~\ref{rem:negative}.
In general,
$\Phi_n$, \smash{$\bar\Phi_n$}, $\Phi_{-n}$, \smash{$\bar\Phi_{-n}$}
are mutually independent except for
$\bar\Phi_{-1} = -\bar\Phi_1$.
The problem is to determine how many of the Voss nets obtained from
these symmetries can be linearly independent, that is, what is
the dimension of the $\RR$-linear space of the Voss nets obtainable by
the Guichard process.

The following examples show how to derive the answer from the
symmetry properties of the initial pseudospherical surface by using
Lemma~\ref{lemma1}.

\begin{Example}
The Beltrami pseudosphere corresponds to the kink
$\phi = 4 \arctan(\mathrm{e}^{x + y})$, the symmetry invariance properties
of which can be expressed by means of the following four easily verifiable
relations
\[
\Phi_2 = \Phi_1, \qquad
\bar\Phi_2 = \bar\Phi_1 - 2 \Phi_1,
\qquad
\Phi_{-1} = \Phi_1, \qquad
\bar\Phi_{-2} = \bar\Phi_1 + 2 \Phi_1.
\]
Applying $\RO$ to the first two and $\RO^{-1}$ to the last two relations,
it follows by induction that all
$\Phi_n$, $\bar\Phi_n$, $\Phi_{-n}$, $\bar\Phi_{-n}$
are linear combinations of $\Phi_1$ and $\bar\Phi_1$.
Consequently, the space of generated symmetries is of dimension two.
Using Lemma~\ref{lemma1}, we easily see that $\Phi_1$ produces a~degenerate Voss net since $\RO \Phi_1 + \Phi_1 = 2 \Phi_1 \in \RO0$.
On the other hand,
$\RO\bar\Phi_1 + \bar\Phi_1 = \bar\Phi_2 + \bar\Phi_1
 = 2\bar\Phi_1 - 2 \Phi_1 \equiv 2\bar\Phi_1$ mod $\RO0$,
so that $\bar\Phi_1$ produces a~nontrivial Voss surface
(namely the right helicoid according to
Proposition~\ref{prop:rhelicoid:R}).
Consequently, the space of generated Voss nets is of dimension
$2 - 1 = 1$.

For the Dini helicoid~\eqref{eq:Dini:phi}, the answer is 1 as well.
\end{Example}

\begin{Example}
The degenerate two-soliton solution
\[
\phi = 4 \arctan\left(\frac{x - y}{\cosh(x + y)}\right),
\]
of the sine-Gordon equation, see equation~\eqref{eq:Kuen:phi},
satisfies
\[
\Phi_3 = 2 \Phi_2 - \Phi_1, \qquad
\bar\Phi_2 = \bar\Phi_1 - 4 \Phi_1,
\qquad
\Phi_{-1} = -\Phi_2 + 2 \Phi_1, \qquad
\bar\Phi_{-2} = \bar\Phi_1 - 4 \Phi_2 + 8 \Phi1.
\]
Applying $\RO$ to the first two and $\RO^{-1}$ to the last two relations,
it follows by induction that all
$\Phi_n$, \smash{$\bar\Phi_n$}, $\Phi_{-n}$,
\smash{$\bar\Phi_{-n}$}, $n \ge 1$,
are linear combinations of $\Phi_1$, $\Phi_2$ and \smash{$\bar\Phi_1$}.
Consequently, the space of generated symmetries is of dimension three.
To answer Guichard's problem with the help of Lemma~\ref{lemma1},
we solve the equation
$(\RO + \mathrm{Id})\bigl(a \Phi_1 + b\Phi_2 + c\bar\Phi_1\bigr) \in \RO0$
for $a,b,c \in \RR$.
The general solution is $a = 0$, $c = -3 b$, which yields a~1-parametric family of symmetries for which the resulting Voss
net turns out to be degenerate.
Thus, the answer is $3 - 1 = 2$ linearly independent Voss nets,
for which we can choose the koru surface
(see Example~\ref{ex:Voss:koru})
and the eared screw (see Example~\ref{ex:Voss:archim}).
\end{Example}

\begin{Example}
The degenerate three-soliton solution
of the sine-Gordon equation is
\[
4 \arctan\left(\frac{(x + y) \cosh(x + y) - (x - y)^2 \sinh(x + y)}{(x - y)^2 + \cosh^2(x + y)}\right) + 4 \arctan(\mathrm{e}^{x + y})
\]
\cite[equation~(2.13)]{JC-FC-AF-1917} and satisfies
\begin{gather*}
\Phi_4 = 3 \Phi_3 - 3 \Phi_2 + \Phi_1, \qquad
\bar\Phi_2 = \bar\Phi_1 - 6 \Phi_1,
\qquad
\Phi_{-1} = -\Phi_3 + 3 \Phi_2 - 3 \Phi_1, \\
\bar\Phi_{-2} = \bar\Phi_1 + 6 \Phi_3 - 18 \Phi_2 + 18 \Phi_1.
\end{gather*}
In analogy to the previous example, $\bar\Phi_1$, $\Phi_3$, $\Phi_2$, $\Phi_1$
form the basis of
the four-dimensional space of generated symmetries.
With the help of Lemma~\ref{lemma1},
we reveal one symmetry for which the resulting Voss net turns out
to be degenerate, namely $\Phi_3 - 4 \Phi_2 + 7\Phi_1$.
Thus, the answer is that the Guichard procedure generates
$4 - 1 = 3$ linearly independent Voss nets.
\end{Example}

The overall length of Guichard's sequences is expected to be $n$ for all
$n$-soliton and, more generally, all $n$-gap solutions~\cite{Z-T-F-1974}.

\section{More inverse operators}
\label{sect:neg}

Once $\RO^{-1}$ yields practically the same results as $\RO$,
the next possibility is to invert $\RO \pm \lambda \mathrm{Id}$,
$\lambda > 0$.
While the operator $\RO$ is invariant with respect to the translations,
it transforms into its multiple $\mathrm{e}^{-2\gamma} \RO$ under the scaling
$x \mapsto \mathrm{e}^\gamma x$, $y \mapsto \mathrm{e}^{-\gamma} y$.
Indeed, the potential $Q$ transforms to~${\mathrm{e}^{-\gamma} Q}$ and
the rest is obvious.

Since $(\RO \pm \lambda\,\mathrm{Id})^{-1}
 = \lambda^{-1} (\lambda^{-1}\RO \pm \mathrm{Id})^{-1}$ and
$\lambda^{-1}\RO$ is related to $\RO$ by the scaling,
as we have just seen,
the inversion of $\RO \pm \lambda\, \mathrm{Id}$ can be reduced to the inversion
of $\RO \pm \mathrm{Id}$ in combination with the scaling transformation.
We shall be interested in $(\RO + \mathrm{Id})^{-1}$.

\begin{Proposition}
\label{prop:sG.iro}
If~$\Psi$ satisfies the equation $\Psi_{xy} = (\cos\phi)\Psi$,
then the system
\begin{gather}
\Phi_{xx} = -\Phi
 + \frac{\phi_x}{\tan\phi} \Phi_x
 + \frac{\phi_x}{\sin\phi} \Phi_y
 + \Psi
 - \frac{\phi_x}{\sin\phi} \Psi_y,
\qquad
\Phi_{xy} = \Phi \cos\phi,
\nonumber
\\
\Phi_{yy} = -\Phi
 + \frac{\phi_y}{\tan\phi} \Phi_y
 + \frac{\phi_y}{\sin\phi} \Phi_x
 - \frac{\phi_y}{\tan\phi} \Psi_y
 + \Psi_{yy}
\label{eq:sG:iro}
\end{gather}
is compatible and
$\Phi \equiv (\RO + \mathrm{Id})^{-1} \Psi$.
\end{Proposition}

\begin{proof}
To recover the operator $(\RO + \mathrm{Id})^{-1}$ in the Guthrie form,
it is enough to express $\Phi$ from the system
\[
\Psi = \Phi_{xx} + \phi_x Q + \Phi,
\qquad
Q_x = \Phi_x \phi_x,
\qquad
Q_y = \Phi \sin \phi.
\]
Eliminating $Q$, we easily arrive at the three
equations~\eqref{eq:sG:iro}.
It is also straightforward to check that the system is compatible
if $\Psi$ satisfies $\Psi_{xy} = \Psi \cos\phi$.
Obviously from the second equation~\eqref{eq:sG:iro},
$\Phi$ is a~symmetry, that is,
the correspondence $\Psi \mapsto \Phi$ is a~recursion operator.
\end{proof}

\begin{Remark}
We can give an alternative proof of Lemma~\ref{lemma1}.
Its proof uses system~\eqref{eq:degen:iR0}, which
coincides with system~\eqref{eq:sG:iro} for $\Psi = 0$,
meaning that
$
\Phi \in (\RO + \mathrm{Id})^{-1} 0
$
by Proposition~\ref{prop:sG.iro}.
The equivalence with equation~\eqref{eq:degen:RO:1} is obvious
by the definition of inversion.
\end{Remark}

\begin{Corollary}
Denoting
\[
\mathbf \Phi = \left(\begin{matrix}
 \Phi \\ \Phi_{x} \\ \Phi_{y}
\end{matrix}\right)
\]
and
\[
\mathbf A~= \left(\begin{matrix}
 0 & 1 & 0 \\
 \llap{$-$}1 & {\phi_x}/{\tan\phi} &
 {\phi_x}/{\sin\phi}
\\
 \cos\phi & 0 & 0
\end{matrix}\right),
\qquad
\mathbf B = \left(\begin{matrix}
 0 & 0 & 1 \\
 \cos\phi & 0 & 0 \\
 \llap{$-$}1 & {\phi_y}/{\sin\phi} &
 \phi_y/{\tan\phi}
\end{matrix}\right),
\]
the recursion operator $(\RO + \mathrm{Id})^{-1}$ can be written in the matrix
form
\begin{equation}
\label{eq:sG:iro:m}
\mathbf \Phi_x
 = \mathbf A~\mathbf \Phi
+ \left(\begin{matrix}
 0 \\
 \Psi - \phi_x\Psi_y/{\sin\phi} \\
 0
\end{matrix}\right),
\qquad
\mathbf \Phi_y
 = \mathbf B \mathbf \Phi
+ \left(\begin{matrix}
 0 \\
 0 \\
 \Psi_{yy} - \phi_y \Psi_y/{\tan\phi}
\end{matrix}\right).
\end{equation}
\end{Corollary}

\begin{proof}
It is straightforward to check that systems~\eqref{eq:sG:iro:m}
and~\eqref{eq:sG:iro} are equivalent.
\end{proof}

\begin{Remark}
\label{rem:lambda}

Relative to the sine-Gordon equation,
$\mathbf A~\dif x + \mathbf B \dif y$
is a~$\mathfrak{gl}_3$-valued zero-curvature representation
independent of $\lambda$,
which turns out to be gauge-equivalent to the adjoint representation
of the well-known $\mathfrak{sl}_2$-valued zero-curvature
representation~\cite{A-K-N-S,Z-T-F-1974} for
$\lambda = 1$.
The full $\lambda$-dependent zero-curvature representation
can be obtained either by applying the scaling symmetry to
$\mathbf A~\dif x + \mathbf B \dif y$ or from the inversion
$(\RO + \lambda\,\mathrm{Id})^{-1}$,
see~\cite[p.~252]{MM-spp} and references therein.
For an interesting connection with Laplace invariants see
\cite{H-Kh-P-2016}.
\end{Remark}

To proceed further, we recall the Gauss--Weingarten system
\begin{gather*}
\ve r_{xx}
 = \frac{\phi_x}{\tan \phi} \ve r_x - \frac{\phi_x}{\sin \phi} \ve r_y,
\qquad
\ve r_{xy} = \sin \phi \ve n,
\qquad
\ve r_{yy}
 = -\frac{\phi_y}{\sin \phi} \ve r_x + \frac{\phi_y}{\tan \phi} \ve r_y,
\\
\ve n_x = \frac{1}{\tan \phi} \ve r_x - \frac{1}{\sin \phi} \ve r_y,
\qquad
\ve n_y = \frac{1}{\tan \phi} \ve r_y - \frac{1}{\sin \phi} \ve r_x
\end{gather*}
for pseudospherical surfaces in Dini's parameterisation $x$, $y$,
corresponding to a~solution $\phi$ of the sine-Gordon equation.

\begin{Proposition} \label{prop:sG:iro}
The recursion operator $(\RO + \mathrm{Id})^{-1}$ can be written in the form
\[
(\RO + \mathrm{Id})^{-1}\Psi = \ve n \cdot \ve q,
\]
where $\ve q$ can be obtained by the quadratures
\begin{gather}
\ve q_{x}
 = \frac{\ve r_y}{\sin\phi}\left(
 \frac{\phi_x}{\sin\phi}\Psi_y - \Psi\right),
\qquad
\ve q_{y}
 = \frac{\ve r_x}{\sin\phi}\left(
 \frac{\phi_y}{\tan\phi} \Psi_y - \Psi_{yy}\right).
\label{eq:sG:iro:q}
\end{gather}
\end{Proposition}

\begin{proof}
Eliminating $\ve r_x$ and $\ve r_y$ from the Gauss--Weingarten
equations~\eqref{eq:sG:GW} via
\[
\ve r_x = -\frac{\ve n_x}{\tan\phi} - \frac{\ve n_y}{\sin\phi},
\qquad
\ve r_y = -\frac{\ve n_x}{\sin\phi} - \frac{\ve n_y}{\tan\phi},
\]
we get
\begin{gather*}
\ve n_{xx} = -\ve n
 + \frac{\phi_x}{\tan \phi} \ve n_x
 + \frac{\phi_x}{\sin \phi} \ve n_y,
\qquad
\ve n_{xy} = \cos \phi\, \ve n,
\\
\ve n_{yy} = -\ve n
 + \frac{\phi_y}{\sin \phi} \ve n_x
 + \frac{\phi_y}{\tan \phi} \ve n_y.
\end{gather*}
In matrix form,
$
\mathbf N_{x}
 = \mathbf A \mathbf N$, $
\mathbf N_{y}
 = \mathbf B
\mathbf N$,
where
\[
\mathbf N = \left(\begin{matrix}
 \ve n \\ \ve n_{x} \\ \ve n_{y}
\end{matrix}\right)
\]
denotes the $3\times3$ matrix composed of the rows
$\ve n$, $\ve n_x$, $\ve n_y$,
the matrices $\mathbf A$, $\mathbf B$ being the same as above.
Substituting
$
\mathbf A = \mathbf N_x \mathbf N^{-1}$, $
\mathbf B = \mathbf N_y \mathbf N^{-1}
$
into equations~\eqref{eq:sG:iro:m}, the recursion operator becomes
\begin{gather*}
\mathbf \Phi_x
 = \mathbf N_x \mathbf N^{-1} \mathbf \Phi
+ \left(\begin{matrix}
 0 \\
 \Psi - \phi_x\Psi_y/{\sin\phi} \\
 0
\end{matrix}\right),
\qquad
\mathbf \Phi_y
 = \mathbf N_x \mathbf N^{-1} \mathbf \Phi
+ \left(\begin{matrix}
 0 \\
 0 \\
 \Psi_{yy} - \phi_y \Psi_y/{\tan\phi}
\end{matrix}\right).
\end{gather*}
Equivalently,
\begin{gather}
\bigl(\mathbf N^{-1} \mathbf \Phi\bigr)_x
 = \mathbf N^{-1} \left(\begin{matrix}
 0 \\
 \Psi - \phi_x\Psi_y/{\sin\phi} \\
 0
\end{matrix}\right),
\qquad
\bigl(\mathbf N^{-1} \mathbf \Phi\bigr)_y
 = \mathbf N^{-1} \left(\begin{matrix}
 0 \\
 0 \\
 \Psi_{yy} - \phi_y \Psi_y/{\tan\phi}
\end{matrix}\right).\label{eq:sG:iro:N}
\end{gather}
Therefore, $\mathbf N^{-1} \mathbf \Phi$ can be identified with the
vector $\ve q$ introduced by system~\eqref{eq:sG:iro:q}.
Now, recall that the $3\times3$ matrix $\mathbf N$ has vectors
$\ve n, \ve n_x, \ve n_y$ as rows.
Then $\mathbf N^{-1}$ is the $3\times3$ matrix
\begin{gather*}
\mathbf N^{-1}
 = \frac{\operatorname{adj} \mathbf N}{\det \mathbf N}
 = \left( \begin{matrix}
 \dfrac{\ve n_x \times \ve n_y}{[\ve n, \ve n_x, \ve n_y]} &
 \dfrac{\ve n_y \times \ve n}{[\ve n, \ve n_x, \ve n_y]} &
 \dfrac{\ve n \times \ve n_x}{[\ve n, \ve n_x, \ve n_y]}\end{matrix}\right)^\top
 = \left( \begin{matrix}\ve n & {-\dfrac{\ve r_y}{\sin\phi}} &
 {-\dfrac{\ve r_x}{\sin\phi}}\end{matrix}\right)^\top
\end{gather*}
column-wise.
Substituting this expression into equation~\eqref{eq:sG:iro:N},
we get the required formulas~\eqref{eq:sG:iro:q} in the first row.
\end{proof}

The following proposition immediately leads to the result of
Example~\ref{ex:nc} above.

\begin{Proposition}
\label{prop:R+I:ideal}
The ideal~$(\RO + \mathrm{Id})^{-1}0$ consists of the projections of the
normal vector $\ve n$ to arbitrary directions, that is, of the dot
products
$
\ve n \cdot \ve c$,
where $\ve c$ is an arbitrary constant vector.
\end{Proposition}

\begin{proof} \hskip 1ex
Substituting $\Psi = 0$ into equation~\eqref{eq:sG:iro:q},
we obtain
$\ve q_x = 0$ and $ \ve q_y = 0$.
Hence the statement.
\end{proof}

Denoting $\tilde\Phi_0 = \ve n \cdot \ve c$, the shadow
$\tilde\Phi_0$ can be interpreted as the starting member
of what we shall call the {\it negative $\RO + \mathrm{Id}$ hierarchy}.

\begin{Proposition} \label{prop:iRO:1}
The first member $(\RO + \mathrm{Id})^{-1}(\ve n \cdot \ve c)$ of the
negative $\RO + \mathrm{Id}$ hierarchy coincides with
$
\tilde\Phi_1 = \ve n \cdot \ve q_1$,
where $\ve q_1$ is determined by
\begin{gather}
\ve q_{1,x}
 = -\frac{\ve r_y}{\sin\phi}
 \left(\ve n - \frac{\phi_x}{\sin\phi} \ve n_y\right) \cdot \ve c,
\qquad
\ve q_{1,y}
 = \frac{\ve r_x}{\sin\phi}
 \left(\ve n - \frac{\phi_y}{\sin\phi} \ve n_x\right) \cdot \ve c.
\label{eq:sG:gw:q:1}
\end{gather}
\end{Proposition}

\begin{proof}
Substituting $\Psi = \ve n \cdot \ve c$ into system~\eqref{eq:sG:iro:q},
we obtain exactly system~\eqref{eq:sG:gw:q:1}.
\end{proof}

\section{Beyond the Guichard sequence}
\label{sect:iRO:Voss}

It turns out that the vector $\ve q$ determined by
formulas~\eqref{eq:sG:iro:q} determines a~Voss surface
away from the Guichard sequence and independent of it.

\begin{Proposition} \label{prop:rtilde}
Let $\ve r(x,y)$ be the position vector of a~pseudospherical surface
in the Dini parameterisation.
Let $\Psi$ be any solution of the Moutard equation
$
\Psi_{xy} = \Psi \cos \phi$.
Let
\[
\tilde X = \frac{\phi_x}{\sin\phi}\Psi_y - \Psi, \qquad
\tilde Y = \frac{\phi_y}{\tan\phi} \Psi_y - \Psi_{yy}.
\]
\begin{itemize}\itemsep=0pt
\item[$(i)$] Introduce a~vector $\tilde{\ve r}$ by
\begin{equation}
\label{eq:sG:iro:rtilde}
\tilde{\ve r}_{x}
 = \frac{\tilde X}{\sin\phi} \ve r_y,
\qquad
\tilde{\ve r}_{y}
 = \frac{\tilde Y}{\sin\phi} \ve r_x.
\end{equation}
Then the net~$\tilde{\ve r}(x,y)$, if nonsingular, is a~Voss net.
The corresponding fundamental forms are
\[
\tilde{\mathrm{I}}
 = \frac{\tilde X^2 \dif x^2
 + 2 \tilde X \tilde Y \cos\phi \dif x \dif y + \tilde Y^2 \dif y^2}
 {\sin^2 \phi},
\qquad
\tilde{\mathrm{II}}
 = \tilde X \dif x^2 + \tilde Y \dif y^2.
\]

\item[$(ii)$] The singular points~$\det \tilde{\mathrm{I}} = 0$
of~$\tilde{\ve r}$ are determined by
$
\tilde X \tilde Y = 0$.

\item[$(iii)$] The functions $\tilde X$, $\tilde Y$ obtained for
$
\Psi = \mathcal R \Phi + \Phi
$
coincide with the functions $X$, $Y$ obtained for $\Phi$
via Proposition~\rm\ref{VsG->Voss}.
\end{itemize}
\end{Proposition}

\begin{proof}
(i) Obviously,
$
\tilde{\ve r}_x \cdot \ve n = \tilde{\ve r}_y \cdot \ve n = 0.
$
Viewing $\tilde{\ve r}$ as the position vector of a~surface,
we see that $\tilde{\ve r}$ and $\ve r$
are related by the parallelism of normals.

Obviously from $\tilde{\mathrm{II}}_{12} = 0$,
the parameters $x$, $y$ are conjugate.
Next, computing
$\tilde{\ve r}_x \times \tilde{\ve r}_{xx}$ and~${\tilde{\ve r}_y \times \tilde{\ve r}_{yy}}$,
one easily sees that they are
parallel to $\ve n_y$ and $\ve n_x$, respectively.
Consequently, both are orthogonal to $\ve n$ and the geodesic curvatures
\[
\tilde{\mathrm{gc}}_1 = \frac{[\tilde{\ve r}_x, \tilde{\ve r}_{xx}, \ve n]}{(\tilde{\ve r}_x \cdot \tilde{\ve r}_x)^{3/2}},
\qquad
\tilde{\mathrm{gc}}_2 = \frac{[\tilde{\ve r}_y, \tilde{\ve r}_{yy}, \ve n]}{(\tilde{\ve r}_y \cdot \tilde{\ve r}_y)^{3/2}}
\]
along the coordinate lines are zero.

(ii) It is easy to check that
$\det \tilde I = \tilde X \tilde Y\mkern-4mu/{\sin^2\phi}$.

(iii)
This can be seen by comparing formulas~\eqref{eq:Voss:Phi:TD}
and~\eqref{eq:sG:iro:rtilde}.
\end{proof}

When applying Proposition~\ref{prop:rtilde} to $\Psi = \ve n \cdot \ve c$,
we obtain $\ve q_1$ of Proposition~\ref{prop:iRO:1}.
The resulting Voss net can be nontrivial and linearly independent of the
nets of the Guichard sequence, i.e., sequence obtained from the pmKdV
hierarchy or the Khor'kova hierarchy via Proposition~\ref{VsG->Voss}.

\begin{Example}[the {\it arch}]
Let us apply Proposition~\ref{prop:iRO:1}
(that is, Proposition~\ref{prop:rtilde} with $\Psi = \ve n \cdot \ve c$)
to the Beltrami pseudosphere,
 that is, to $\ve r$ given by equation~\eqref{eq:PS:r}
and $\phi$ given by equation~\eqref{eq:PS:phi}.
Let $\ve c \in \mathbb S^2$ be an arbitrary unit vector.
Choosing $\ve c = [0,0,1]$, we get the constant
$\tilde{\ve r} = \ve 0$.
Choosing $\ve c = [\sin\alpha, \cos\alpha, 0]$ orthogonal to $[0,0,1]$,
we obtain a~family of Voss nets rotated around the axis $[0,0,1]$
by $\alpha$.
Thus, we are left with a~single Voss net
\begin{gather*}
\tilde{\ve r}
 = \left[\frac{1 - \cos(2 x - 2 y)}{4 \tanh(x + y)} - \frac{x + y}{2},
 - \frac{\sin(2 x - 2 y)}{4 \tanh(x + y)},
 \frac{1}{2} \cos(x - y) \cosh(x + y)\right]
\end{gather*}
up to the rigid motions.
Computing the loci of singular points according to part~(ii) of
Proposition~\ref{prop:rtilde},
we get the curves
$
\tanh(x + y) = \pm \tan (x - y)$,
known as the pseudo-elliptic radioids~\cite{ET-1939}.
The part of the Voss surface delimited by two pseudo-elliptic radioids,
$
\tanh(x + y) < |{\tan}(x - y)|$,
is shown in Figure~\ref{fig:arch}.
The ``spinal curve'' $y = x$ is the catenary.
Therefore, this part of the surface can be viewed as a~bent lath
version of the common catenary arch.
It would be interesting to investigate to what extent it can
support its own weight without collapsing.

\begin{figure}\centering
\includegraphics[scale=0.5]{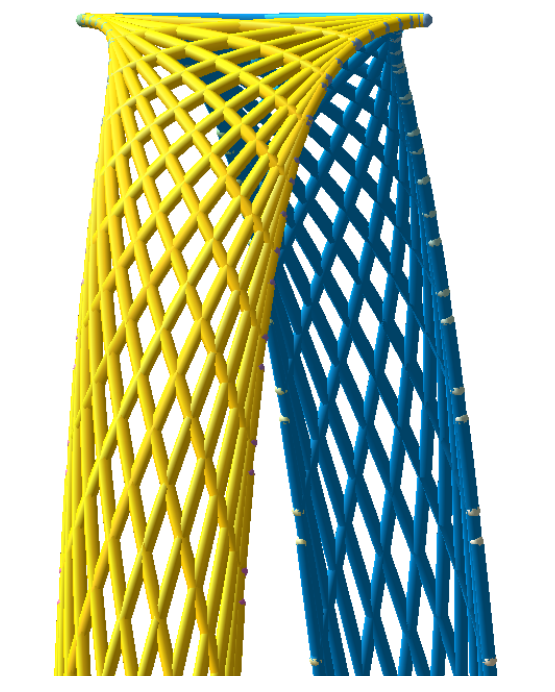} \qquad
\includegraphics[scale=0.5]{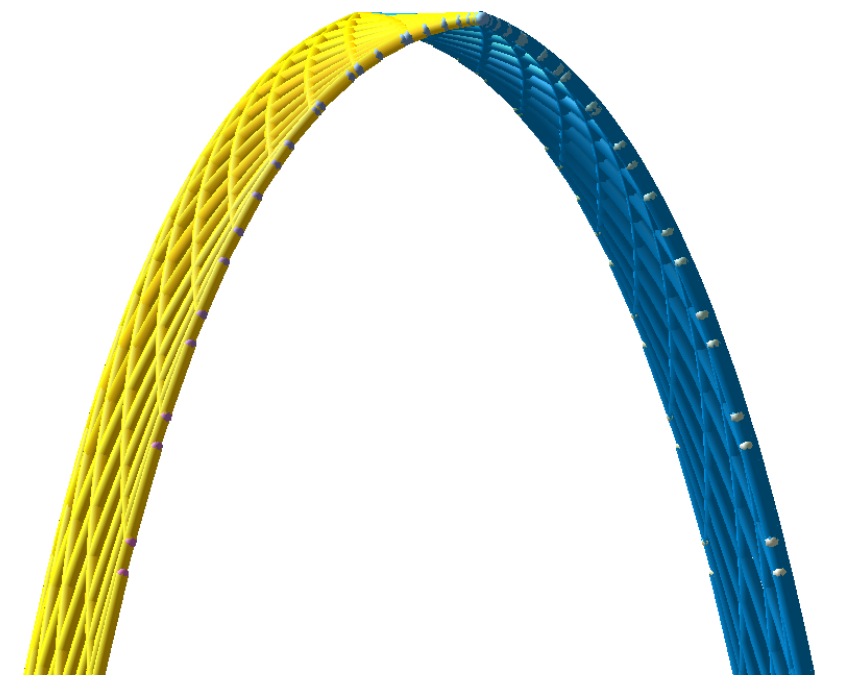}
\caption{The Voss arch from two different viewpoints.}
\label{fig:arch}
\end{figure}
\end{Example}

We finish this paper with a~result analogous to Lemma~\ref{lemma1}
in the context of Proposition~\ref{prop:rtilde}.

\begin{Lemma} \label{lemma2}
Formula~\eqref{eq:sG:iro:rtilde}
turns a~solution $\Psi$ of the
Moutard equation into a~degenerate Voss surface
if and only if
$
\Psi \in \RO 0$.
\end{Lemma}

\begin{proof}
The lemma is a~consequence of Lemma~\ref{lemma1} and
Proposition~\ref{prop:rtilde}\,(iii).
 A~short direct proof is as follows.
By Proposition~\ref{prop:rtilde}\,(ii),
the surface~\eqref{eq:sG:iro:rtilde}
is degenerate if and only if $\tilde X = 0$ and $\tilde Y = 0$.
The overdetermined system formed by
$\tilde X = 0$, $\tilde Y = 0$ and the Moutard equation for~$\Psi$
reduces to the passive form
\[
\Psi_{x} = \frac{\phi_{xx}}{\phi_x} \Psi,\qquad
\Psi_{y} = \frac{\phi_{xy}}{\phi_x} \Psi,\]
with the general solution
$\Psi = \text{const}\, \phi_x$ or, equivalently,
$\Psi \in \RO 0$.
\end{proof}

\section{The division algebra of recursion operators}
\label{sect:dig}

Recursion operators form an associative $\RR$-algebra
with respect to operations defined by
\[
(\U \circ \V)\Phi = \U\V\Phi, \qquad
(\U + \V)\Phi = \U \Phi + \V \Phi, \qquad
(r \U)\Phi = r \U\Phi, \qquad r \in \RR.
\]
Previously, this algebra was studied for the heat equation
\cite{SDK-ROP-2023,ISK-AMV-1989,JK-AV-RV-2017}
(is non-commutative),
the Krichever--Novikov and the Landau--Lifshitz equation~\cite{D-S-2008}.
For evolution equations like these, ambiguities do not pose a~serious problem.
To cope with ambiguities discussed in Section~\ref{sect:sG:RO},
we factor out the ideal formed by the trivial operators in the sense of
Definition~\ref{def:RO:triv}.

\begin{Proposition}
Trivial recursion operators form a~two-sided ideal in the associative
$\RR$-algebra of all recursion operators.
\end{Proposition}

\begin{proof}
Indeed, if $\Z$ is trivial, then $\Z\U\Phi = \Z\Phi$ by definition.
Moreover, if
$\Z\Phi = W^{(1)} \Sigma^{(1)} + \cdots + W^{(n)} \Sigma^{(n)}$, then
$\U\Z\Phi = W^{(1)} \U\Sigma^{(1)} + \cdots + W^{(n)} \U\Sigma^{(n)}$
by linearity.
\end{proof}

The two-sided ideal induces a~congruence relation on the $\RR$-algebra
of all recursion operators.

\begin{Definition}
Two recursion operators $\U$ a~$\V$ are said to be {\it congruent},
if $\U - \V$ is trivial.
We write $\U \equiv \V$.
\end{Definition}

Now, consider the associative algebra $\langle\RO\rangle$ generated by
the operator~\eqref{eq:sG:RO:Guthrie}.
Let $\mathfrak R$ denote the quotient algebra
$\langle\RO\rangle/{\equiv}$.
Elements of $\mathfrak R$ are recursion operators modulo the congruence
$\equiv$.

In consequence of Proposition~\ref{prop:sG:RO:Guthrie:inv},
we get the following corollary.

\begin{Corollary}
As elements of $\mathfrak R$, the operators $\RO$, $\RO'$ are
mutually inverse.
\end{Corollary}

\begin{Proposition}
The associative algebra $\mathfrak R$ is commutative and
isomorphic to the algebra $\mathbb R[\RO]$ of polynomials in $\RO$.
\end{Proposition}

\begin{proof}
Obviously, the associative algebra
$(\mathfrak R, +, 0, -, \circ, \mathbb R)$
generated by the single element $\RO$ is commutative
and consists of polynomials
$\mathcal P
 = \RO^n + \nu_{n-1} \RO^{n-1} + \cdots + \nu_1 \RO + \nu_0 \,\mathrm{Id}$.
To show that $\mathfrak R \cong \RR[\RO]$, we have to verify
that no non-constant polynomial $\mathcal P$ is congruent to zero.
To this end, consider $\mathcal P \phi_x = \RO^n \phi_x
 + \nu_{n-1} \RO^{n-1} \phi_x + \cdots + \nu_1 \RO \phi_x$.
According to Remark~\ref{rem:negative}, symmetries
$\RO^{i}\phi_x$ are members of the pmKdV
hierarchy of ever increasing order.
Consequently, $\mathcal P \phi_x$ is nonzero in~$\mathfrak R$.
Hence the statement.
\end{proof}

Finally, consider the division algebra $\mathfrak Q$ of $\mathfrak R$.
By the previous proposition, $\mathfrak Q$ is isomorphic to~$\RR(\RO)$,
the division algebra of $\RR[\RO]$.
The elements of $\RR(\RO)$
can be viewed as formal fractions~${\P/\Q}$,
where $\P,\Q \in \RR[\RO]$, $\Q \not\equiv 0$.
Since not all recursion operators are invertible
(for instance, the trivial ones),
some fractions need not correspond to true recursion operators,
in principle.
The next proposition proves that they actually do.

\begin{Proposition}
All operators $\P \circ \Q^{-1}$,
corresponding to formal fractions $\P/\Q$, where $\Q \not\equiv 0$,
are B\"acklund auto-transformations of the Moutard equation.
\end{Proposition}

\begin{proof}
It suffices to show that all inverses $(\RO + \lambda\,\mathrm{Id})^{-1}$,
corresponding to linear factors $\RO + \lambda\,\mathrm{Id}$ of $\Q$,
are true B\"acklund transformations.
For $\lambda = 0$,
this follows from Proposition~\ref{prop:sG:RO:Guthrie:inv}.
For~${\lambda = 1}$, see Proposition~\ref{prop:sG.iro},
which easily extends to arbitrary $\lambda \ne 0$.
\end{proof}

The entire division algebra will be employed in a~future work.

\subsection*{Acknowledgements}

Special thanks are due to Ivan Izmestiev for illuminating discussions
and for generously providing a~copy of his work in progress. A~discussion with Adam Doliwa was no less important.

Thanks are also due to the organizers of the conferences
{\it Differential Geometry and its Applications}, Brno, Sept. 8--12, 2025,
and
{\it 18th Symposium on Integrable Systems}, Warsaw, July 6--7, 2026,
where the main results were presented.

Last but not least, the author extends his gratitude to the anonymous
referees for their time, valuable feedback, and insistence on rigour.

The author was supported by the RVO institutional funding under
I\v{C} 47813059.


\pdfbookmark[1]{References}{ref}
\LastPageEnding


\begin{thebibliography}{19}
\footnotesize\itemsep=0pt
\bibitem{A-K-N-S}
Ablowitz M.J., Kaup D.J., Newell A.C., Segur H., Method for solving the
 sine-{G}ordon equation,
 \href{https://doi.org/10.1103/PhysRevLett.30.1262}{\textit{Phys. Rev. Lett.}}
 \textbf{30} (1973), 1262--1264.

\bibitem{Ai-1983}
Aiyer R.N., Recursion operators for infinitesimal transformations and their
 inverses for certain nonlinear evolution equations,
 \href{https://doi.org/10.1088/0305-4470/16/2/008}{\textit{J.~Phys.~A}}
 \textbf{16} (1983), 255--262.

\bibitem{Bab-1996}
Babich M.V., Methods of the theory of sine-{G}ordon equation and the {V}oss
 surfaces, {P}reprint NTZ 46/1996, Universit\"at Leipzig, Zentrum f\"ur
 H\"ohere Studien, 1996,
 \url{https://cds.cern.ch/record/320244/files/SCAN-9702082.pdf}.

\bibitem{Bia-I}
Bianchi L., Lezioni di geometria differenziale. {V}ol.~I, E.~Spoerri, Pisa,
 1894.

\bibitem{Bia-II}
Bianchi L., Lezioni di geometria differenziale. {V}ol.~II, E.~Spoerri, Pisa,
 1903.

\bibitem{B-V-V}
Bocharov A.V., Chetverikov V.N., Duzhin S.V., Khorkova N.G., Krasilshchik I.S.,
 Samokhin A.V., Torkhov~Yu.N., Verbovetsky A.M., Vinogradov A.M., Symmetries
 and conservation laws for differential equations of mathematical physics,
 \textit{Transl. Math. Monogr.}, Vol. 182,
 \href{https://doi.org/10.1090/mmono/182}{American Mathematical Society},
 Providence, RI, 1999.

\bibitem{JC-FC-AF-1917}
Cen J., Correa F., Fring A., Degenerate multi-solitons in the sine-{G}ordon
 equation,
 \href{https://doi.org/10.48550/arXiv.1705.04749}{\textit{J.~Phys.~A}}
 \textbf{50} (2017), 435201, 20~pages,
 \href{http://arxiv.org/abs/nlin.SI/1705.04749}{arXiv:nlin.SI/1705.04749}.

\bibitem{D-S-2008}
Demskoi D.K., Sokolov V.V., On recursion operators for elliptic models,
 \href{https://doi.org/10.48550/arXiv.nlin/0607071}{\textit{Nonlinearity}}
 \textbf{21} (2008), 1253--1264,
 \href{http://arxiv.org/abs/nlin.SI/0607071}{arXiv:nlin.SI/0607071}.

\bibitem{JBD-1913}
Diller J.B., \"Uber die den {E}nneperschen {F}l\"achen konstanten negativen
 {K}rummungsmasses entsprechenden {V}o{\ss}schen {F}l\"achen, Ph.D.~Thesis,
 {U}niversit\"at W\"urzburg, 1913.

\bibitem{UD-1870}
Dini U., Sopra alcune formole generali della teoria delle superficie, e loro
 applicazioni, \href{https://doi.org/10.1007/BF02420031}{\textit{Ann. Mat.
 Pura Appl.}} \textbf{4} (1870), 175--206.

\bibitem{D-B-1977}
Dodd R.K., Bullough R.K., Polynomial conserved densities for the sine-{G}ordon
 equations, \href{https://doi.org/10.1098/rspa.1977.0012}{\textit{Proc. Roy.
 Soc. London Ser.~A}} \textbf{352} (1977), 481--503.

\bibitem{AD-1999}
Doliwa A., Integrable discrete geometry with ruler and compass, in Symmetries
 and {I}ntegrability of {D}ifference {E}quations ({C}anterbury, 1996),
 \textit{London Math. Soc. Lecture Note Ser.}, Vol.~255,
 \href{https://doi.org/10.1017/CBO9780511569432.011}{Cambridge University
 Press}, Cambridge, 1999, 122--136.

\bibitem{LPE-1906}
Eisenhart L.P., Associate surfaces,
 \href{https://doi.org/10.1007/BF01449817}{\textit{Math. Ann.}} \textbf{62}
 (1906), 504--538.

\bibitem{LPE-1909}
Eisenhart L.P., A~treatise on the differential geometry of curves and surfaces,
 Ginn, Boston, 1909.

\bibitem{LPE-1914}
Eisenhart L.P., Transformations of surfaces of~{V}oss,
 \href{https://doi.org/10.2307/1988587}{\textit{Trans. Amer. Math. Soc.}}
 \textbf{15} (1914), 245--265.

\bibitem{LPE-1915}
Eisenhart L.P., Conjugate systems with equal point invariants,
 \href{https://doi.org/10.2307/2007145}{\textit{Ann. of Math.}} \textbf{18}
 (1916), 7--17.

\bibitem{LPE-1917}
Eisenhart L.P., Certain surfaces of {V}oss and surfaces associated with them,
 \href{https://doi.org/10.1007/BF03014896}{\textit{Rend. Circ. Mat. Palermo}}
 \textbf{42} (1917), 145--166.

\bibitem{SF-1939}
Finikoff S., D\'eformation \`a r\'eseau conjugu\'e persistant et probl\`emes
 g\'eom\'etriques qui s'y rattachent, \textit{M\'emor. Sci. Math.}, Vol.~96,
 \href{https://www.numdam.org/item/MSM_1939__96__1_0/}{Gauthier}, Paris, 1939.

\bibitem{BG-1926}
Gambier B., D\'eformation continue d'un h\'elico\"{\i}de en h\'elico\"{\i}de
 avec r\'eseau conjugu\'e permanent, \textit{Bull. Sci. Math.} \textbf{50}
 (1926), 308--328.

\bibitem{BG-1927}
Gambier B., Surfaces de {V}oss et {G}uichard: {S}urfaces associ\'ees et
 adjointes. {D}\'eformation avec r\'eseau conjugu\'e permanent,
 \href{https://doi.org/10.1007/BF02545661}{\textit{Acta Math.}} \textbf{51}
 (1928), 83--131.

\bibitem{BG-1931}
Gambier B., Surfaces de~{V}oss--{G}uichard,
 \href{https://doi.org/10.24033/asens.813}{\textit{Ann. Sci. \'{E}cole Norm.
 Sup.}} \textbf{48} (1931), 359--396.

\bibitem{CG-1890b}
Guichard C., Recherches sur les surfaces \`a courbure totale constante et sur
 certaines surfaces qui s'y rattachent, \href{https://www.numdam.org/item/10.24033/asens.333.pdf}{\textit{Ann. Sci. \'Ecole Norm. Sup.}} \textbf{7} (1890), 233--264.

\bibitem{CG-1890a}
Guichard C., Sur une classe particuli\`ere d'\'equations aux d\'eriv\'ees
 partielles dont les invariants sont \'egaux, \href{https://www.numdam.org/item/10.24033/asens.333.pdf}{\textit{Ann. Sci. \'Ecole Norm. Sup.}} \textbf{7} (1890), 19--22.

\bibitem{GAG-1993}
Guthrie G.A., Miura transformations and symmetry groups of differential
 equations, Ph.D.~Thesis, {U}niversity Canterbury, 1993,
 \url{https://ir.canterbury.ac.nz/server/api/core/bitstreams/43340926-fe51-4a80-9113-2bded2ddd046/content}.

\bibitem{GAG-1994}
Guthrie G.A., Recursion operators and non-local symmetries,
 \href{https://doi.org/10.1098/rspa.1994.0094}{\textit{Proc. Roy. Soc. London
 Ser.~A}} \textbf{446} (1994), 107--114.

\bibitem{H-Kh-P-2016}
Habibullin I.T., Khakimova A.R., Poptsova M.N., On a~method for constructing
 the {L}ax pairs for nonlinear integrable equations,
 \href{https://doi.org/10.1088/1751-8113/49/3/035202}{\textit{J.~Phys.~A}}
 \textbf{49} (2016), 035202, 35~pages,
 \href{http://arxiv.org/abs/1506.02563}{arXiv:1506.02563}.

\bibitem{HKF}
Hundertwasser F., Hundertwasser koru flag, hundertwasser Art Centre,
 Whang\=arei, New Zealand,
 \url{https://www.hundertwasserartcentre.co.nz/about/hundertwasser/hundertwasser-koru-flag/}.

\bibitem{Izm}
Izmestiev I., Raffaeli M., Rasoulzadeh A., Smooth, discrete, and semi-discrete
 {V}oss surfaces, {i}n preparation.

\bibitem{NGK-1988}
Khor'kova N.G., Conservation laws and nonlocal symmetries,
 \href{https://doi.org/10.1007/BF01158125}{\textit{Math. Notes}} \textbf{44}
 (1988), 562--568.

\bibitem{KNRRS-2024}
Kilian M., Nawratil G., Raffaelli M., Rasoulzadeh A., Sharifmoghaddam K.,
 Interactive design of discrete {V}oss nets and simulation of their rigid
 foldings, \href{https://doi.org/10.1016/j.cagd.2024.102346}{\textit{Comput.
 Aided Geom. Design}} \textbf{111} (2024), 102346, 24~pages.

\bibitem{K-K-2002}
Konno K., Kakuhata H., A~fused hierarchy,
 \href{https://doi.org/10.1016/S0034-4877(02)80022-2}{\textit{Rep. Math.
 Phys.}} \textbf{49} (2002), 239--248.

\bibitem{AVK-2020}
Kostin A.V., On~{D}ini helicoids in the~{M}inkowski space,
 \href{https://doi.org/10.1007/s10958-023-06772-9}{\textit{J.~Math. Sci.}}
 \textbf{276} (2023), 517--524.

\bibitem{SDK-ROP-2023}
Koval S.D., Popovych R.O., Point and generalized symmetries of the heat
 equation revisited,
 \href{https://doi.org/10.1016/j.jmaa.2023.127430}{\textit{J.~Math. Anal.
 Appl.}} \textbf{527} (2023), 127430, 21~pages,
 \href{http://arxiv.org/abs/2208.11073}{arXiv:2208.11073}.

\bibitem{ISK-PHMK-2000}
Krasil'shchik I.S., Kersten P.H.M., Symmetries and recursion operators for
 classical and supersymmetric differential equations, \textit{Math. Appl.},
 Vol.~507, \href{https://doi.org/10.1007/978-94-017-3196-6}{Springer}, Dordrecht, 2000.

\bibitem{ISK-AMV-1989}
Krasil'shchik I.S., Vinogradov A.M., Nonlocal trends in the geometry of
 differential equations: symmetries, conservation laws, and {B}\"acklund
 transformations, \href{https://doi.org/10.1007/BF00131935}{\textit{Acta Appl.
 Math.}} \textbf{15} (1989), 161--209.

\bibitem{JK-AV-RV-2017}
Krasil'shchik J., Verbovetsky A., Vitolo R., The symbolic computation of
 integrability structures for partial differential equations, \textit{Texts Monogr.
 Symbol. Comput.}, \href{https://doi.org/10.1007/978-3-319-71655-8}{Springer},
 Cham, 2017.

\bibitem{RL-RC-2018}
Lin R., Conte R., On a~surface isolated by Gambier,
 \href{https://doi.org/10.48550/arXiv.1805.10450}{\textit{J.~Nonlinear Math.
 Phys.}} \textbf{25} (2018), 509--514,
 \href{http://arxiv.org/abs/1805.10450}{arXiv:1805.10450}.

\bibitem{L-C-1993}
Lou S.Y., Chen W.Z., Inverse recursion operator of the {AKNS} hierarchy,
 \href{https://doi.org/10.1016/0375-9601(93)90677-R}{\textit{Phys. Lett.~A}}
 \textbf{179} (1993), 271--274.

\bibitem{L-H-L-2021}
Lou S.Y., Hu X.B., Liu Q.P., Duality of positive and negative integrable
 hierarchies via relativistically invariant fields,
 \href{https://doi.org/10.1007/jhep07(2021)058}{\textit{J.~High Energy Phys.}}
 \textbf{2021} (2021), no.~7, 058, 21~pages,
 \href{http://arxiv.org/abs/2103.04549}{arXiv:2103.04549}.

\bibitem{MM-alro}
Marvan M., Another look on recursion operators, in Differential {G}eometry and
 {A}pplications, Proc. Conf. Brno., 393--402, \url{https://emis.muni.cz/proceedings/6ICDGA/IV/index.html}.

\bibitem{MM-rzcr}
Marvan M., Reducibility of zero curvature representations with application to
 recursion operators,
 \href{https://doi.org/10.1023/B:ACAP.0000035588.67805.0b}{\textit{Acta Appl.
 Math.}} \textbf{83} (2004), 39--68,
 \href{http://arxiv.org/abs/nlin.SI/0306006}{arXiv:nlin.SI/0306006}.

\bibitem{MM-spp}
Marvan M., On the spectral parameter problem,
 \href{https://doi.org/10.1007/s10440-009-9450-4}{\textit{Acta Appl. Math.}}
 \textbf{109} (2010), 239--255,
 \href{http://arxiv.org/abs/0804.2031}{arXiv:0804.2031}.

\bibitem{MDTFB-2020}
Montagne N., Douthe C., Tellier X., Fivet C., Baverel O., Voss surfaces: A~design space for geodesic gridshells, in Inspiring the {N}ext {G}eneration
 (Proc. IASS Annual Symp. 2020 and 7th Int. Conf. Spatial Structures),
 \href{https://doi.org/10.15126/900337}{Spatial Structures Research Centre of
 the University of Surrey}, 2021, 3473--3483.

\bibitem{MDTFB-2022}
Montagne N., Douthe C., Tellier X., Fivet C., Baverel O., Discrete {V}oss
 surfaces: {D}esigning geodesic gridshells with planar cladding panels,
 \href{https://doi.org/10.1016/j.autcon.2022.104200}{\textit{Autom. Constr.}}
 \textbf{140} (2022), 104200, 19~pages.

\bibitem{PJO-1977}
Olver P.J., Evolution equations possessing infinitely many symmetries,
 \href{https://doi.org/10.1063/1.523393}{\textit{J.~Math. Phys.}} \textbf{18}
 (1977), 1212--1215.

\bibitem{CJP-1991}
Papachristou C.J., Lax pair, hidden symmetries, and infinite sequences of
 conserved currents for self-dual {Y}ang--{M}ills fields,
 \href{https://doi.org/10.1088/0305-4470/24/17/015}{\textit{J.~Phys.~A}}
 \textbf{24} (1991), L1051--L1055.

\bibitem{P-M-2007}
Popov A.G., Maevskii E.V., Analytical approaches to the study of the
 sine-{G}ordon equation and pseudospherical surfaces,
 \href{https://doi.org/10.1007/s10958-007-0183-5}{\textit{J.~Math. Sci}}
 \textbf{142} (2007), 2377--2418.

\bibitem{AR-2025}
Rasoulzadeh A., Classification of smooth alignable {V}oss surfaces,
 \href{http://arxiv.org/abs/2603.17501}{arXiv:2603.17501}.

\bibitem{AR-1889}
Razzaboni A., Delle superficie sulle quali due serie di geodetiche formano un
 sistema conjugato, \textit{Mem. Acad. Sci. Bologne} \textbf{9} (1889),
 765--776.

\bibitem{R-L-W-2001}
Reid G.J., Lin P., Wittkopf A.D., Differential elimination-completion
 algorithms for {DAE} and {PDAE},
 \href{https://doi.org/10.1111/1467-9590.00159}{\textit{Stud. Appl. Math.}}
 \textbf{106} (2001), 1--45.

\bibitem{sammlung}
Roth-Kleyer G. (Ed.), G\"ottingen collection of mathematical models and instruments,
 \url{http://modellsammlung.uni-goettingen.de/}.

\bibitem{VR-2000}
Rovenski V., Geometry of curves and surfaces with {MAPLE},
 \href{https://doi.org/10.1007/978-1-4612-2128-9}{Birkh\"auser}, Boston, MA,
 2000.

\bibitem{ES-1924}
Salkowski E., Zur {T}heorie der {V}o\ss schen und der {G}uichardschen
 {F}l\"achen, \href{https://doi.org/10.1007/BF01188086}{\textit{Math.~Z.}}
 \textbf{20} (1924), 241--252.

\bibitem{S-W-2001b}
Sanders J.A., Wang J.P., Integrable systems and their recursion operators,
 \href{https://doi.org/10.1016/S0362-546X(01)00630-7}{\textit{Nonlinear
 Anal.}} \textbf{47} (2001), 5213--5240.

\bibitem{S-W-2001a}
Sanders J.A., Wang J.P., On recursion operators,
 \href{https://doi.org/10.1016/S0167-2789(00)00188-3}{\textit{Phys.~D}}
 \textbf{149} (2001), 1--10.

\bibitem{RS-1970}
Sauer R., Differenzengeometrie,
 \href{https://doi.org/10.1007/978-3-642-86411-7}{Springer}, Berlin, 1970.

\bibitem{RS-HG-1931}
Sauer R., Graf H., \"{U}ber {F}l\"achenverbiegung in {A}nalogie zur
 {V}erknickung offener {F}acettenflache,
 \href{https://doi.org/10.1007/BF01455828}{\textit{Math. Ann.}} \textbf{105}
 (1931), 499--535.

\bibitem{AT-1903}
Tachauer A., \"Uber diejenigen {R}otationsfl\"achen, auf denen zwei {S}charen
 geod\"atischer {L}inien ein konjugiertes {S}ystem bilden, \textit{Arch. Math.
 Phys.} \textbf{6} (1903), 60--84.

\bibitem{ET-1939}
Turri\`ere E., Une nouvelle courbe de transition pour les raccordements
 progressifs: la radio\"{\i}de pseudoelliptique,
 \href{https://doi.org/10.24033/bsmf.1291}{\textit{Bull. Soc. Math. France}}
 \textbf{67} (1939), 62--99.

\bibitem{JMV-1991}
Verosky J.M., Negative powers of {O}lver recursion operators,
 \href{https://doi.org/10.1063/1.529234}{\textit{J.~Math. Phys.}} \textbf{32}
 (1991), 1733--1736.

\bibitem{AV-1888}
Voss A., Diejenigen fl\"achen, auf denen zwei schaaren geod\"atischer linien
 ein conjugirtes system bilden, \textit{Sitzungsber. Bayer. Akad. Wiss.
 Math.-Naturw. Klasse} \textbf{18} (1888), 95--102.

\bibitem{Z-1994}
Zadadaev S.A., Solutions of traveling waves type to the sine-{G}ordon equation,
 and pseudodifferential surfaces, \textit{Moscow Univ. Math. Bull.}
 \textbf{49} (1994), 34--40.

\bibitem{VEZ-BGK-1984}
Zakharov V.E., Konopelchenko B.G., On the theory of recursion operator,
 \href{https://doi.org/10.1007/BF01403883}{\textit{Comm. Math. Phys.}}
 \textbf{94} (1984), 483--509.

\bibitem{Z-T-F-1974}
Zakharov V.E., Takhtajan L.A., Faddeev L.D., A~complete description of the
 solutions of the ``sine-{G}ordon'' equation, \textit{Sov. Phys. Dokl.}
 \textbf{19} (1975), 824--826.

\end{thebibliography}
\end{document}